\documentstyle[12pt]{article}
 \input amssym.def
\input amssym
\begin{document}
\def \Z{\Bbb Z}
\def \C{\Bbb C}
\def \R{\Bbb R}
\def \Q{\Bbb Q}
\def \N{\Bbb N}
\def \wt{{\rm wt}}
\def \tr{{\rm tr}}
\def \span{{\rm span}}
\def \Res{{\rm Res}}
\def \Res{{\rm QRes}}
\def \End{{\rm End}}
\def \E{{\rm End}}
\def \Ind {{\rm Ind}}
\def \Irr {{\rm Irr}}
\def \Aut{{\rm Aut}}
\def \Hom{{\rm Hom}}
\def \mod{{\rm mod}}
\def \ann{{\rm Ann}}
\def \<{\langle}
\def \>{\rangle}
\def \t{\tau }
\def \a{\alpha }
\def \e{\epsilon }
\def \l{\lambda }
\def \L{\Lambda }
\def \g{\gamma}
\def \b{\beta }
\def \om{\omega }
\def \o{\omega }
\def \c{\chi}
\def \ch{\chi}
\def \cg{\chi_g}
\def \ag{\alpha_g}
\def \ah{\alpha_h}
\def \ph{\psi_h}
\def \be{\begin{equation}\label}
\def \ee{\end{equation}}
\def \bl{\begin{lem}\label}
\def \el{\end{lem}}
\def \bt{\begin{thm}\label}
\def \et{\end{thm}}
\def \bp{\begin{prop}\label}
\def \ep{\end{prop}}
\def \br{\begin{rem}\label}
\def \er{\end{rem}}
\def \bc{\begin{coro}\label}
\def \ec{\end{coro}}
\def \bd{\begin{de}\label}
\def \ed{\end{de}}
\def \pf{{\bf Proof. }}
\def \voa{{vertex operator algebra}}

\newtheorem{thm}{Theorem}[section]
\newtheorem{prop}[thm]{Proposition}
\newtheorem{coro}[thm]{Corollary}
\newtheorem{conj}[thm]{Conjecture}
\newtheorem{lem}[thm]{Lemma}
\newtheorem{rem}[thm]{Remark}
\newtheorem{de}[thm]{Definition}
\newtheorem{hy}[thm]{Hypothesis}
\makeatletter
\@addtoreset{equation}{section}
\def\theequation{\thesection.\arabic{equation}}
\makeatother
\makeatletter

\baselineskip=24pt
\begin{center}{\Large \bf Vertex operator algebras associated to
admissible representations of  $\hat{sl}_{2}$}

\vspace{0.5cm}
Chongying Dong\footnote{Supported by NSF grant DMS-9303374 and a
research grant from the Committee on Research, UC Santa Cruz.},
Haisheng Li and Geoffrey Mason\footnote{Supported by NSF grant
DMS-9401272 and a research grant from the Committee on Research, UC
Santa Cruz.}\\ Department of Mathematics, University of California,
Santa Cruz, CA 95064
\end{center}

\section{Introduction}
Let $\{e,f,h\}$ be a standard basis for $\frak g=sl_{2}$ such that
$[e,f]=h$, $[h,e]=2e$, $[h,f]=-2f$, $\hat\frak g$ the corresponding
affine Lie algebra and $L(\ell,j)$ the irreducible highest weight
$\hat\frak g$-module of level $\ell$ with highest weight $j$.  It is
well known that the vacuum representation $L(\ell,0)$ has a natural
vertex operator algebra (or chiral algebra) structure for any $\ell\ne
-2$ (cf. [FZ]). If $\ell$ is a positive integer, the chiral algebra
$L(\ell,0)$ of the WZNW models in the content of conformal field
theory has been well understood. For example, the fusion rules are
obtained by using primary field decomposition (cf. [GW], [TK]) or by
Verlinde formula (cf. [K], [V]), $n$ point functions are calculated
[KZ].

In the content of vertex operator algebra, it has been proved
(cf. [DL], [FL], [Li1]) that any $\Z_+$-graded weak $L(\ell,0)$-module
is completely reducible and the set of equivalence classes of
irreducible $L(\ell,0)$-modules is the set of equivalence classes of
standard $\hat\frak g$-modules of level $\ell$.  Thus $L(\ell, 0)$ is
rational (defined in Section 2).  (It has been proved recently in
[DLiM2] that {\em any weak} $L(\ell,0)$-module is completely reducible
and the set of equivalence classes of irreducible {\em weak}
$L(\ell,0)$-modules is the set of equivalence classes of standard
$\hat\frak g$-modules of level $\ell$.)  The modular invariance of the
vector space linearly spanned by the characters $tr_{L(\ell,j)}e^{2\pi
i\tau(L(0)-{c_{\ell}\over 24})}$ for all standard modules of level
$\ell$ is obtained in [KP] by using the explicit character formulas or
follows from a general theorem of Zhu [Z].  The fusion rules
are computed in [FZ] by studying certain associative algebras and its
bimodules associated to $L(\ell,0)$ and its irreducible modules.

If $\ell$ is rational such that $\displaystyle{\ell+2={p\over q}}$ for
some coprime positive integers $p\ge 2$ and $q$, Kac and Wakimoto
[KW1]-[KW2] found finitely many distinguished irreducible representations,
called admissible (or modular invariant) representations.
In this case the fusion rules among admissible modules have
been calculated in the content of conformal field theory (cf. [AY],
[BF], [MW]) by employing different methods, but different methods
sometime give different results. Especially, Verlinde formula gives
negative fusion rules.

If $j$ is not an integer, $\tr_{L(\ell,j)}e^{2\pi
i\tau(L(0)-{c_{\ell}\over 24})}$ does not exist (because the
homogeneous subspaces are infinite-dimensional) so that the character
$\tr_{L(\ell,j)}e^{2\pi i\tau(L(0)-{1\over 2}zh(0)-{c_{\ell}\over
24})}$ was considered in [KW1]-[KW2], where $z$ is a positive rational
number less than $1$.  In [KW1], a formula for $\tr_{L(\ell,j)}e^{2\pi
i\tau(L(0)-{1\over 2}zh(0)-{c_{\ell}\over 24})}$ in terms of theta
functions was given and a transformation law under
$S(\tau,z)=(-\tau^{-1}, -z\tau)$ were found.  Later, the
transformation law was corrected by adding an extra factor
[KW2]. After this correction, it is not clear that the space linearly
spanned by all $\tr_{L(\ell,j)}e^{2\pi i\tau(L(0)-{1\over
2}zh(0)-{c_{\ell}\over 24})}$ for admissible weights $j$ is invariant
under the action of the modular group $PSL(2,{\Z})$, where
$S(\tau,z)=(-\tau^{-1}, -z\tau)$.

Our modest purpose of this paper is to study these admissible
representations from the point of view of vertex operator algebra.
We show that $L(\ell,0)$ is a $\Q$-graded rational vertex operator
algebra under a new Virasoro algebra and its irreducible modules
are exactly these admissible modules for $\hat\frak g.$ We extend
Zhu's $A(V)$-theory ([FZ],[Z]) to $\Q$-graded vertex operator algebras
and apply this theory to $L(\ell,0)$ to calculate the fusion rules.
The new Virasoro algebra also gives a natural interpretation of
the characters $\tr_{L(\ell,j)}e^{2\pi i\tau(L(0)-{1\over
2}zh(0)-{c_{\ell}\over 24})}.$

We explain these results in details in the following.
In the first part of the paper, we prove that all
the admissible representations of level $\ell$ constitute the set of
irreducible $\Z_+$-graded weak $L(\ell,0)$-modules among all the highest-weight
irreducible $\hat\frak g$-modules of level $\ell$.
 This has been
implicit in references such as [AY], [BF], [FM].
It follows from this
result and a complete reducibility theorem of Kac-Wakimoto [KW2] that
any weak $L(\ell,0)$-module from category $\cal{O}$ is completely
reducible. Let $N_{+}$ be the sum of all positive root spaces of
$\hat{\frak{g}}$. Let $\cal{E}$ be the category of weak
$L(\ell,0)$-modules $W$ on which $N_{+}$ is locally nilpotent, {\it
i.e.}, for any $u\in W,$ there is a positive integer $k$ such that
$N_{+}^{k}u=0$.  Then we prove that any weak $L(\ell,0)$-module from
category $\cal{E}$ is completely reducible.

Since some admissible weights are not integers if $\ell$ is not
integral, Zhu's algebra $A(L(l,0))$ [Z] has infinite-dimensional
irreducible modules. This implies that $L(\ell,0)$ is not rational and
that Zhu's $C_{2}$-condition (another crucial condition for Zhu's
theorem of modular invariance) is not true either.

In the second part of the paper, we study
the vertex operator algebra $L(\ell,0)$ under a new Virasoro
algebra. Let $\omega$ be the original Segal-Sugawara Virasoro vector
of $L(\ell,0)$. Set $\omega_{z}=\omega +{1\over 2}zh(-2){\bf 1}\in
L(\ell,0)$, where $z$ is a complex number.  Then $\omega_{z}$ is a
Virasoro vector of $L(\ell,0)$ with a central charge
$c_{\ell,z}=c_{\ell}-6\ell z^{2}$ and
$L_{z}(0)=(\omega_{z})_{1}=L(0)-{1\over 2}zh(0)$. If we choose $z=0,
{1\over 2}$, we obtain the homogeneous grading and the rescaled
principal grading (cf. [K], [LW]), respectively.  Let $z$ be positive
rational number less than $1$.  Note that the vertex operator algebra
$(L(\ell,0),Y, {\bf 1},\omega_{z})$ is ${\Q}$-graded instead of
${\Z}$-graded. We extend Zhu's $A(V)$-theory of one-to-one
correspondence [Z] between the set of equivalence classes of
irreducible admissible $V$-modules and the set of equivalence classes
of irreducible $A(V)$-modules and Frenkel-Zhu's $A(M)$ theory [FZ] for
fusion rules to any ${\Q}$-graded vertex operator algebra. It follows
from our complete reducibility theorem in the first part that any
${\Q}_{+}$-graded weak $L(\ell,0)$-module under the new Virasoro
vector $\omega_{z}$ is completely reducible.  That is,
$(L(\ell,0),\omega_{z})$ is rational.  By using the
Malikov-Feigin-Fuchs's singular vector expressions [MFF] and the
Fuchs' projection formula [F] we find all the fusion rules and prove
that $(L(\ell,0),\omega_{z})$ satisfies the $C_{2}$-finite
condition. Our results on fusion rule agree with the corresponding
results in [BF].

It is natural for us to consider the modified characters $\tr e^{2\pi
i\tau (L_{z}(0)-{1\over 24}c_{\ell,z})}$, that is, $\tr e^{2\pi i\tau
(L(0)-{1\over 2}zh(0)-{1\over 24}(c_{\ell}-6\ell z^{2}))}$.  Using
KW's character formula [KW1] we find that these modified characters
are modular functions so that $c_{\ell,z}$ is the modular anomaly
rather than $c_{\ell}$.  (Then the characters $\tr e^{2\pi i\tau
(L(0)-{1\over 2}zh(0)-{1\over 24}c_{\ell})}$ are obviously not modular
functions.)  One may ask: Is the space linearly spanned by our new
characters invariant under the transformation $S (\tau)=-\tau^{-1}$
with $z$ being fixed?  This will be discussed in our coming paper [DLiM3].

We should mention that the vertex operator algebras
associated to irreducible highest weight representations of certain
rational levels
for affine Lie algebra $C_n^{(1)}$ have been studied in [A].

\section{Vertex operator algebras $L(\ell,0)$ associated to $\hat{sl}_{2}$}

A vertex operator algebra, or briefly a VOA, is a ${\Z}$-graded
vector space $V=\oplus_{n\in {\Z}}V_{n}$ satisfying a number of axioms.
We refer the reader to [B], [FLM] and [FHL]
for the details of the definition.  However, we would like
to give the definitions of weak modules
and modules in details. A {\em weak $V$-module} is a
pair $(W,Y_{W})$, where $W$ is a vector space and $Y_{W}(\cdot,z)$ is a linear
map from $V$ to
$({\rm End} W)[[z,z^{-1}]]$
satisfying the following axioms:
(1) $Y_{W}({\bf 1},z)=id_{W}$; (2) $Y_{W}(a,z)u\in W((z))$ for any $a\in V,u\in
W$;
(3) $Y_{W}(L(-1)a,z)={d\over dz}Y_{W}(a,z)$ for $a\in V$; (4) the Jacobi
identity:
\begin{eqnarray}
&
&z_{0}^{-1}\delta\left(\frac{z_{1}-z_{2}}{z_{0}}\right)Y_{M}(a,z_{1})Y_{M}(b,z_{2})u
-z_{0}^{-1}\delta\left(\frac{z_{2}-z_{1}}{-z_{0}}\right)Y_{M}(b,z_{2})Y_{M}(a,z_{1})u
\nonumber\\
&
&=z_{2}^{-1}\delta\left(\frac{z_{1}-z_{0}}{z_{2}}\right)Y_{M}(Y(a,z_{0})b,z_{2})u
\end{eqnarray}
for $a,b\in V, u\in W$. A weak $V$-module $(W,Y_{W})$ is called a {\em
$V$-module} if $L(0)$ semisimply acts on $W$ with the decomposition
into $L(0)$-eigenspaces $M=\oplus_{h\in {\C}}M_{h}$ such that for any
$h\in {\C}$, $\dim M_{h}<\infty, M_{h+n}=0$ for $n\in {\Z}$
sufficiently small.

A {\em ${\Z}_{+}$-graded weak} $V$-module [FZ] is a weak $V$-module $W$
together with
a $\Z_+$-gradation $W=\oplus_{n=0}^{\infty}W(n)$ such that
\begin{eqnarray}
a_{m}W(n)\subseteq W(k+n-m-1)\;\;\;\mbox{ for }a\in V_{k}, m,n\in {\Z},
\end{eqnarray}
where $W(n)=0$ by definition for $n<0$. One may define the notions of
``submodule''and ``irreducible submodule'' accordingly.  A VOA $V$ is
said to be {\em rational} if any ${\Z}_{+}$-graded weak $V$-module is
a direct sum of irreducible ${\Z}_{+}$-graded weak $V$-modules. It was
proved in [DLiM1] that if $V$ is rational, there are only finitely many
irreducible ${\Z}_{+}$-graded weak $V$-modules up to equivalence and
any irreducible weak $V$-module is a module.

Let $\{e,f,h\}$ be the standard basis for $\frak g=sl_{2}$ with the
 commutation relations:
$[e,f]=h. [h,e]=2e, [h,f]=-2f$. We fix the normalized Killing form on
$\frak g$ such that $\<h,h\>=2$.
Let $\tilde{\frak g}=\tilde{sl}_{2}={\C}[x,x^{-1}]\otimes \frak g+{\C}c$
be the affine Lie algebra and identify $\frak g$ with $x^0\otimes\frak g.$
Set $a(n)=a\otimes x^n$ for $a\in \frak g$ and $n\in\Z$ for convenience.
Define subalgebras
\begin{eqnarray}
& & N_{+}={\C}e+x{\C}[x]\otimes \frak g,\;\;
N_{-}={\C}f+x^{-1}{\C}[x^{-1}]\otimes \frak g,\label{d2.3}\\
& & B=N_+\oplus \C h\oplus \C c,\ \ \ \ P=\C[x]\otimes\frak g\oplus \C c.
\end{eqnarray}
Then $\tilde{sl}_{2}=N_{+}\oplus {\C}h\oplus {\C}c\oplus N_{-}$.

Let $\hat{\frak{g}}=\tilde{\frak{g}}\oplus \C d$ be the extended affine algebra
[K],
where
$$[d,c]=0,\; [d,x^{n}\otimes a]=n(x^{n}\otimes a)\;\;\;\mbox{for }a\in
\frak{g}, n\in {\Z}.$$
Let $H=\C h\oplus \C c\oplus \C d$ be the Cartan subalgebra of
$\hat{\frak{g}}$,
$\alpha_{0},\alpha_{1}$ be the simple roots of $\hat{\frak{g}}$, let
$\Gamma_{+}=\Z_{+}\alpha_{0}\oplus \Z_{+}\alpha_{1}$ and let
$\Lambda_{0},\Lambda_{1}$
be the fundamental weights of $\hat{\frak{g}}$ [K].
Let $\bar{\rho}$ be half of the sum of positive roots of $\frak{g}$
and set $\rho=\bar{\rho}+2\Lambda_{0}$ [K].
For any $\lambda\in H^{*}$, denote by $M(\lambda)$ (resp.  $L(\lambda)$) the
Verma (resp. the
irreducible highest weight) $\hat{\frak{g}}$-module. When restricted to
$\tilde{\frak{g}}$,
$L(\lambda)$ is
an irreducible $\tilde{\frak{g}}$-module [K]. It is clear that $L(\lambda)$ and
$L(\mu)$ are
isomorphic $\tilde{\frak{g}}$-module if and only if $\lambda\in
\mu+{\C}\delta$.
 As commonly
used in many references, we use the notation $L(\ell, j)$ for the
$\tilde{\frak{g}}$-module
$L(\lambda)$, where
$\ell=\<\lambda,c\>, j=\<\lambda,\alpha_{1}\>=j$.
Conversely, let $M$ be a restricted $\tilde{\frak{g}}$-module of level $\ell\ne
-2$. Then
we extend $M$ to a $\hat{\frak{g}}$-module
by letting $d$ act on $M$ as $-L(0)$. In this paper we shall consider any
restricted $\tilde{\frak{g}}$-module as a $\hat{\frak{g}}$-module in this way.

For a complex number $l$ and a $\C h$-module $U$ which can be regarded as
a $B$-module by $N_+$ acting trivially and $c$ acting as $l,$
let $M(\ell,U)$ be the {\em generalized Verma} $\tilde{\frak g}$-{\em module}
$U({\tilde{\frak g}})\otimes_{U(B)}U$ [Le] of level $l$ or Weyl module.
If $U=\C$ is one-dimensional $\C h$-module on which
$h$ acts as a fixed complex number $j$ the corresponding
module is an ordinary {\em Verma module} denoted by $M(l,j).$
Note that $U$ can be identify with
the subspace $1\otimes_{U(B)}U$ of $M(l,U).$ Then  $M(l,j)$ has
a unique maximal submodule which intersects trivially with $\C$ and
$L(\ell,j)$  is isomorphic to the corresponding irreducible highest weight
module.

Similarly, one can define the
generalized Verma $\tilde{\frak g}$-module $V(\ell,U)=U(\tilde{\frak
g})\otimes_{U(P)}U$ for any $\frak g$-module $U$ which can be extended
to a $P$-module by setting $x\C[x]\otimes \frak g$ acting trivially
and $c$ acting as $l.$ Note that if $U=\C$ is the trivial $\frak g$-module
then $V(\ell,\C)$ is a quotient of $M(l,0)$ and $L(l,0)$ is the irreducible
quotient of $V(\ell,\C)$ modulo the unique maximal submodule which intersects
$\C$ trivially.

It is well-known that $V(\ell,{\C})$ and $L(\ell,0)$
have natural vertex operator algebra structures for any $\ell\ne -2$
and that any $M(\ell,U)$ is a weak module for vertex operator algebra
$V(\ell,{\C})$ (cf. [FZ] and [Li1]).

We recall the following Kac-Kazhdan reducibility criterion [KK]:

\begin{prop}\label{p2.5} The Verma module
$M(\ell,j)$ is reducible if and only if there are some positive integers $n, k$
such that
 one of the three conditions hold:
\begin{eqnarray}\label{e2.17}
(I)\;\; j=n-1-(k-1)t;\;\;\;(II)\;\; j=-n+kt;\;\;\;(III)\;\; \ell+2=0,
\end{eqnarray}
where $t=\ell+2$.
\end{prop}

\begin{rem}\label{r2.6}  Since any restricted $\tilde{\frak g}$-module of level
$\ell$ is a weak $V(\ell,{\C})$-module ([FZ], [Li1]), $V(\ell,{\C})$ is always
irrational.
If $t=\ell +2\not\in \Q_+$, then it follows from
Proposition \ref{p2.5} that
$V(\ell,{\C})=L(\ell,0).$ Therefore, $L(\ell,0)$ is an
irrational vertex operator algebra.
\end{rem}

Recall from [KW1] that a weight $\lambda\in H^{*}$ is said to be {\em
admissible} if the following conditions hold:

(i) $\<\lambda+\rho,\alpha\> >0$ for all but finitely many positive roots
$\alpha$ of $\hat{\frak{g}}$;

(ii) $\<\lambda+\rho,\alpha\> \notin \{-1,-2,\cdots\}$ for any positive root
$\alpha$ of
$\hat{\frak{g}}$;

(iii) The set of positive roots $\alpha$ satisfying $\<\lambda+\rho,\alpha\>\in
{\Z}_{+}$ spans a
$2$-dimensional subspace of $H^{*}$.

A complex number $\ell$ is called {\em an admissible level} if there is an
admissible weight $\lambda$
such that $\<\lambda,c\>=\ell$.
It was proved in [KW1] that $\ell$ is an admissible level if and only if
$\ell=-2+{p\over q}$,
where $p$ and $q$ are coprime positive integers with $p\ge 2$ and $j$ is an
admissible weight of
level $\ell$ if and only if
$j=n-k{p\over q}$ for some $n, k\in {\Z}_{+}, n\le p-2, k\le q-1.$
{}From now on we will assume that $t=\ell + 2={p\over q}$, where $p$
and $q$ are coprime positive integers with $p\ge 2$.

\begin{rem}\label{r2.7}
Let $j=n-kt$ be an admissible weight. Then $j=-(p-n)+(q-k)t$. Since $p$ and $q$
are
relatively prime, $r-st=0$ for $r,s\in {\Z}, 0\le s\le q-1$
if and only if $s=0, r=0$.
Consequently,
the expression $j=n-kt$ of an admissible weight $j$ with
$n,k\in {\Z}_{+}, n\le q-2, k\le q-1$ is unique.
\end{rem}

A vector $w$ in a highest weight module $M$ for $\tilde\frak g$
is called a {\em singular vector}
 if $w$ is a highest weight vector which generates a proper submodule.
It is well known that the singular vectors of $M(\ell,j)$
give the key information
for determining the module structure of $L(l,j)$ and the fusion rules.
In [MFF] an expression for singular vectors in terms of non-integral
powers of elements of $\hat\frak g$
was found as follows (see [MFF]
for details):

\begin{prop}\label{pmff} {\rm [MFF]} Let $j=n-1-(k-1)t$ where $n$ and $k$ are
positive
integers
satisfying $1\le n\le p-1, 1\le k\le q$ and let $v$ be a highest weight vector
of the
Verma module $M(\ell,j)$. Set
\begin{eqnarray}\label{es}
&
&F_{1}(n,k)=f(0)^{n+(k-1)t}e(-1)^{n+(k-2)t}f(0)^{n+(k-3)t}e(-1)^{n+(k-4)t}\nonumber\\
& &\hspace{2 cm}\cdots e(-1)^{n-(k-2)t}f(0)^{n-(k-1)t},\\
&
&F_{2}(n,k)=e(-1)^{p-n+(q-k)t}f(0)^{p-n+(q-k-1)t}e(-1)^{p-n+(q-k-2)t}f(0)^{p-n+(q-k-3)t}\nonumber\\
& &\hspace{2 cm}\cdots f(0)^{p-n-(q-k+1)t}e(-1)^{p-n-(q-k)t}.
\end{eqnarray}
Then $v_{j,1}=F_{1}(n,k)v, v_{j,2}=F_{2}(n,k)v$ are singular
vectors of $M(\ell,j)$
of degrees $n((k-1)\alpha_{0}+k\alpha_{1})$ and
$(p-n)((q+1-k)\alpha_{0}+(q-k)\alpha_{1})$, respectively. Moreover, the
maximal proper submodule of $M(\ell,j)$ is generated by $v_{j,1}$ and
$v_{j,2}.$
\end{prop}

\begin{rem}\label{r2.5} Note that $v_{0,2}=F_2(1,1){\bf 1}$ generates
the maximal proper submodule of $V(\ell,\C).$
\end{rem}

For any complex number $\alpha$, following [F] and [FM] we set
$H_{\alpha}=fe-\alpha h- \alpha (\alpha +1)$. Then
\begin{eqnarray}
& &H_{\alpha}H_{\beta}=H_{\beta}H_{\alpha},\;\;\; e^{m}H_{\alpha}=H_{\alpha -m}
e^{m},\;\;\;f^{m}H_{\alpha}=H_{\alpha +m}f^{m},\label{2.7}\\
& &f^{m}e^{m}=H_{0}H_{1}\cdots H_{m-1},\;
e^{m}f^{m}=H_{-1}H_{-2}\cdots H_{-m},\label{2.8}\\
& &h^{m}e^{n}=e^{n}(h+2n)^{m},\;\; h^{m}f^{n}=f^{n}(h-2n)^{m}
\end{eqnarray}
for any complex numbers $\alpha, \beta$ and for any positive integers $m,n$.

Let $\sigma$ be the anti-automorphism of $U(\frak g)$ such that
$\sigma (a)=-a$ for any $a\in \frak g$. Then $\sigma
(H_{\alpha})=H_{-(\alpha +1)}$ for any complex number $\alpha$.
Let $P_{1}$ be the projection $\tilde{\frak g}$ onto $\frak g$ such that
$P_{1}(t^{n}\otimes a)=a$ for any $a\in \frak g$ and $P_1(c)=0.$

\begin{prop}\label{pf1} [F] The following projection formulas hold:
\begin{eqnarray}
& &P_{1}(F_{1}(n,k))
= \left(\prod_{r=0}^{n-1}\prod_{s=1}^{k-1}H_{r+st}\right)f^{n},\\
& &P_{1}(F_{2}(n,k))
=\left(\prod_{r=1}^{p-n}\prod_{s=1}^{q-k}H_{-r-st}\right)e^{p-n}.
\end{eqnarray}
\end{prop}

Let $B_{0}=\C(f(-1)+f(0))+{\C}[x^{-1}](x^{-2}+x^{-1})\otimes \frak g$. Then
$B_{0}$ is an ideal of
$N_{-}$
such that $N_{-}/B_{0}={\C}T_{+}+{\C}T_{0}+{\C}T_{-},$ denoted by $L_0,$
where $T_{+}=e(-1)+B_{0},
T_{0}=h(-1)+B_{0}, T_{-}=f+B_{0}$, satisfies the following commutation
relations:
\begin{eqnarray}\label{d2.13}
[T_{+}, T_{-}]=T_{0},\;\; [T_{0},T_{+}]=-2T_{+},\;\; [T_{0},T_{-}]=2T_{-}.
\end{eqnarray}
Let $P$ be the natural quotient map from $U(N_{-})$ onto $U(L_{0})$. For any
complex number $\alpha$,
we  define $G_{\alpha}=T_{-}T_{+}-\alpha T_{0}+\alpha (\alpha+1)$. Then
\begin{eqnarray}
& &G_{\alpha}G_{\beta}=G_{\beta}G_{\alpha},\;\;\; T_{+}^{m}G_{\alpha}=G_{\alpha
-m}
T_{+}^{m},\;\;\;T_{-}^{m}G_{\alpha}=G_{\alpha +m}T_{-}^{m},\label{d2.14}\\
& &T_{-}^{m}T_{+}^{m}=G_{0}G_{1}\cdots G_{m-1},\;\;\;
T_{+}^{m}T_{-}^{m}=G_{-1}G_{-2}\cdots G_{-m}\label{d2.15}
\end{eqnarray}
for any complex numbers $\alpha, \beta$ and for any positive integer $m$.
Using the same method as suggested in [F] we obtain

\begin{prop}\label{pf3} The following formulas hold:
\begin{eqnarray}
& &P(F_{1}(n,k))
= \left(\prod_{r=0}^{n-1}\prod_{s=1}^{k-1}G_{r+st}\right)T_{-}^{n},\\
& &P(F_{2}(n,k))
=\left(\prod_{r=1}^{p-n}\prod_{s=1}^{q-k}G_{-r-st}\right)T_{+}^{p-n}.
\end{eqnarray}
\end{prop}

Recall that $N_{-}={\C}f+x^{-1}{\C}[x^{-1}]\otimes \frak g$. Set
$B_{2}={\C}f(-1)+x^{-2}{\C}[x^{-1}]\otimes \frak g$. The it is clear that
$B_{2}$ is an
ideal of $N_{-}$. Let $L_{2}=N_{-}/B_{2}$ be the quotient Lie algebra. Then
$L_{2}$ is a
three-dimensional Heisenberg Lie algebra with relations:
$[\bar{e},\bar{f}]=\bar{h},
[\bar{h},\bar{e}]=[\bar{h},\bar{f}]=0$, where $\bar{e}=e(-1)+B_{2},
\bar{f}=f+B_{2}, \bar{h}=h(-1)+B_{2}$.
Let $P_{2}$ be the natural quotient map from $U(N_{-})$ to $U(L_{2})$. Then

\begin{prop}\label{pf2} [F] For any positive integers $1\le n\le p-1,
1\le k \le q$, we have
\begin{eqnarray}
& &P_{2}(F_{1}(n,k))
= \left(\prod_{r=0}^{n-1}\prod_{s=1}^{k-1}\bar{H}_{r+st}\right)\bar{f}^{n},\\
& &P_{2}(F_{2}(n,k))
=\left(\prod_{r=1}^{p-n}\prod_{s=1}^{q-k}\bar{H}_{-r-st}\right)\bar{e}^{p-n},
\end{eqnarray}
where $\bar{H}_{\alpha}=\bar{e}\bar{f}-\alpha\bar{h}$.
\end{prop}

For a $\C h$-module define a linear functional on
$U^{*}\otimes M(\ell,U)$ as follows:
\begin{eqnarray}\label{2.19}
\langle u', u\rangle =u'({\cal{P}}(u))\;\;\;\mbox{ for }u'\in U^{*}, u\in
M(\ell,U),
\end{eqnarray}
where ${\cal{P}}$ is  the projection of $M(\ell,U)$ onto the
subspace $U$. Define
\begin{eqnarray}
I=\{ u\in M(\ell,U)|\langle u', xu\rangle =0 \;\;\;\mbox{ for any
}u'\in U^{*}, x\in U(\tilde{\frak g})\}.
\end{eqnarray}
 It is clear that
$I$ is the unique maximal submodule which
intersects with $U$ trivially. Set $L(\ell,U)=M(\ell,U)/I$ and regard $U$
as a subspace in a natural way. Then $\cal P$ induces a projection
of $L(l,U)$ to $U,$ which is still be denoted by $\cal P,$  and the formula
(\ref{2.19}) also define a linear functional on $U^*\otimes L(l,U).$
Then (see [FZ] or [Li2]) $M(\ell,U)$ and $L(\ell,U)$ are
weak modules for vertex operator algebra $V(\ell,{\C})$. Let $Y(\cdot,z)$
be the vertex operators defining the module structure on $L(\ell,{\C})$.
It is clear that $Y(\cdot,z)$ is an intertwining
operator of type $\left(\!\begin{array}{c}L(\ell,U)\\
V(\ell,{\C})\,L(\ell,U)\end{array}\!
\right)$
(see [FHL] for the definition of intertwining operator).
Let $\cal{Y}$$(\cdot,z)$ be
the intertwining operator of type
$\left(\!\begin{array}{c}L(\ell,U)\\L(\ell,U)\,V(\ell,{\C})\end{array}\!\right)$ defined by
$\cal{Y}$$(u,z)v=e^{zL(-1)}Y(v,-z)u$ (cf. [FHL]).

\begin{lem}\label{l2.8} The $\tilde{\frak g}$-module $L(\ell,U)$ is a
weak module for the vertex operator algebra $L(\ell,0)$ if
and only if
\begin{eqnarray}\label{2.21}
\langle u',{\cal{Y}} (u,z)v_{0,2}\rangle=0\;\;\;\mbox{ {\it for
any }}u'\in U^{*}, u\in U(\frak g)\subset L(\ell,U).
\end{eqnarray}
\end{lem}

{\bf Proof.} It is clear that the condition is necessary. Now we assume
that (\ref{2.21}) holds.
Let $J$ be
the maximal submodule of $V(\ell,\C)$ which intersects $\C$ trivially.
Then $J=U(N_{-})v_{0,2}$. From the
definition of the bilinear form we get
\begin{eqnarray}\label{2.22}
\< u', a{\cal{Y}} (u,z)w\>=0\;\;\;\mbox{ for }u'\in U^{*}, u\in L(\ell,U),
a\in N_{-}U(N_{-}), w\in L(\ell,0).
\end{eqnarray}
By using  the commutator formula
\begin{equation}\label{2.22'}
[a(m), {\cal Y}(u,z)]=\sum_{j\geq 0}{m\choose j}{\cal Y}(a(j)u,z)z^{m-j}
\end{equation}
for $a\in \frak g,$ $m\in \Z$ and $u\in L(l,U)$
together with (\ref{2.21}) we get
\begin{eqnarray}\label{2.23}
\langle u',{\cal{Y}}(u,z)J\rangle=0\;\;\;\mbox{ for any }u'\in U^{*},
u\in U(\frak g)U\subset L(\ell,U).
\end{eqnarray}
{}From the Jacobi identity for the vertex operators against the intertwining
operator we have
\begin{equation}\label{2.24}
{\cal Y}(a(n)u,z)=\sum_{j\geq 0}{n\choose j}a(n-j){\cal Y}(u,z)z^j
-(-1)^n\sum_{j\geq 0}{n\choose j}{\cal Y}(u,z)a(j)z^{n-j}
\end{equation}
for $u\in L(l,U), a\in \frak g$ and $n\in \Z$.
Note that $L(l,U)$ is generated by $U$ as $\tilde{ \frak g}$-module. Combining
(\ref{2.22}), (\ref{2.23}) and (\ref{2.24}) gives
\begin{eqnarray}\label{2.25}
\langle u',{\cal{Y}}(u,z)J\rangle=0\;\;\;\mbox{ for any }u'\in U^{*},
u\in L(\ell,U).
\end{eqnarray}
By the commutator formula (\ref{2.22'}) again, we obtain
\begin{eqnarray}
\langle u',x{\cal{Y}} (u,z)J\rangle=0\;\;\;\mbox{ for
any }u'\in U', x\in U(\tilde{\frak g}), u\in L(\ell,U).
\end{eqnarray}
By the definition of $L(\ell,U)$ we have $\cal{Y}$$(v,z)u=0$ for any $v\in J,
u\in L(\ell,U)$.
Thus, $\cal{Y}$$(\cdot,z)$ induces an intertwining operator of type
 $\left(\begin{array}{c}L(\ell,U)\\L(\ell,U)\,L(\ell,0)\end{array}
\right)$. This proves that $L(\ell,0)$ is a weak module for
$L(\ell,0)$.$\;\;\;\;\Box$

\begin{prop}\label{p2.8}
The $L(\ell,U)$
is a weak $L(\ell,0)$-module if and only if
$f(h)U=0$, where
$$f(h)=\displaystyle{\prod_{r=0}^{p-2}\prod_{s=0}^{q-1}(h-r+st)}.$$
\end{prop}

{\bf Proof.} Recall that ${\cal{Y}}(\cdot,z)$ is the
corresponding nonzero intertwining operator of type
$\left(\!\begin{array}{c}L(\ell,U)\\L(\ell,U)\,V(\ell,\C)\end{array}\!\right)$.
For $n\in {\Z}, a\in \frak g$ we define $\deg (x^{n}\otimes a)=n$. By
(\ref{2.22'}) we obtain
\begin{eqnarray}\label{2.25'}
\langle u',{\cal{Y}}(u,z)av\rangle=\langle u',z^{\deg a}{\cal{Y}}(\sigma
P_{1}(a)u,z)v\rangle
\end{eqnarray}
for $u'\in U^{*}, u\in U(\frak g)U\subseteq L(\ell,U), a\in U(N_{-}),
v\in L(\ell,0)$. Let $a=F_{2}(1,1)$. Then $v_{0,2}=a{\bf 1}.$ By Lemma
\ref{l2.8} and (\ref{2.25'}) $L(l,U)$ is a weak $L(l,0)$-module if and only if
\begin{eqnarray}
\langle u',{\cal{Y}}(\sigma P_{1}(a)u,z){\bf 1}\rangle =0.
\end{eqnarray}
By Proposition \ref{pf1}, we have
\begin{eqnarray}
P_{1}(a)=\prod_{r=1}^{p-1}\prod_{s=1}^{q-1}H_{-r-st}e^{p-1}.
\end{eqnarray}
Then from (\ref{2.7})
\begin{eqnarray}
\sigma P_{1}(x)=(-1)^{p-1}e^{p-1}\prod_{r=1}^{p-1}\prod_{s=1}^{q-1}H_{r-1+st}
=(-1)^{p-1}\prod_{r=1}^{p-1}\prod_{s=1}^{q-1}H_{-p+r+st}e^{p-1}.
\end{eqnarray}
Note that
\begin{eqnarray*}
& &\ \ \ \langle u',{\cal{Y}}(\sigma P_{1}(a)u,z){\bf 1}\rangle\\
& &=\langle u',e^{zL(-1)}(\sigma P_{1}(a))u\rangle\\
& &=\langle u',(\sigma P_{1}(a))u\rangle.
\end{eqnarray*}
Thus $L(\ell,U)$ is a weak $L(\ell,0)$-module if and only if
\begin{eqnarray}
\langle
u',\prod_{r=1}^{p-1}\prod_{s=1}^{q-1}H_{-p+r+st}e^{p-1}U(g)U\rangle=0\;\;\;\mbox{for any }
u'\in U^{*}.
\end{eqnarray}
{}From the grading restriction on the bilinear pair, the later is
equivalent to
$$\displaystyle{\left(\prod_{r=1}^{p-1}\prod_{s=1}^{q-1}H_{-p+r+st}e^{p-1}f^{p-1}
\right)U=0}.$$ By (\ref{2.8}) and the fact that $eU=0$ we have
\begin{eqnarray}
& &\ \ \ \
\prod_{r=1}^{p-1}\prod_{s=1}^{q-1}\prod_{i=1}^{p-1}H_{-p+r+st}H_{-i}U\nonumber\\
& &=\prod_{r=1}^{p-1}\prod_{s=1}^{q-1}\prod_{i=1}^{p-1}(p-r-st)(h-p+r+1+st)i
(h-i+1)U=0.
\end{eqnarray}
Since $p-r-st\ne 0$ for any $1\le r\le p-1, 1\le s\le q-1$ (from Remark
\ref{r2.7}), we obtain
\begin{eqnarray}
\prod_{r=1}^{p-1}\prod_{s=1}^{q-1}\prod_{i=1}^{p-1}(h-p+1+r+st)(h-i+1)U=0.
\end{eqnarray}
Thus
\begin{eqnarray}
\prod_{r=0}^{p-2}\prod_{s=0}^{q-1}(h-r+st)U=0.
\end{eqnarray}
This finishes the proof. $\;\;\;\;\Box$

\begin{coro}\label{c2.8}
The highest weight $\hat{sl}_{2}$-module $L(\ell,j)$
is a weak $L(\ell,0)$-module  if and only if
$j=r-st$ for $0\le r\le p-2,0\le s\le q-1$.
That is, $\hat\frak g$-module $L(\ell,j)$
is a weak $L(\ell,0)$-module  if and only if
$j$ is admissible.
\end{coro}

Let $j$ be an admissible weight so that $L(\ell,j)$ is a weak
$L(\ell,0)$-module. It follows from [FHL] that $L(\ell,0)'$ is a weak
$L(\ell,0)$-module. But $(L(\ell,j)')'\ne L(\ell,j)$ because
$L(\ell,j)$ has infinite-dimensional homogeneous subspaces in
general.  By using the well-known principal grading (cf. [K]), any
$L(\ell,j)=\oplus_{m,n\in {\Z}}L(\ell,j)_{(m,n)}$ becomes a
${\Z}^{2}$-graded space such that each homogeneous subspace is
finite-dimensional. Let $L(\ell,j)^{c}=\oplus_{(m,n)\in {\Z}\times
{\Z}}L(\ell,j)_{m,n}^{*}$ be the restricted dual of $L(\ell,j)$ with
respect to this ${\Z}^{2}$-grading. Then it is clear that
$L(\ell,j)^{c}$ is an irreducible weak $L(\ell,0)$-module satisfying
$(L(\ell,j)^{c})^{c}=L(\ell,j)$. But the lowest $L(0)$-weight subspace
of $L(\ell,j)^{c}$ is a lowest weight $\frak g$-module with $-j$ as
its lowest weight.  Then there is a non-trivial intertwining operator
of type
$\left(\begin{array}{c}L(\ell,0)\\L(\ell,j)L(\ell,j)^{c}\end{array}\right)$
so that $L(\ell,j)$ and $L(\ell,j)^{c}$ are conjugate each other from
the physical point of view.

\begin{coro} Let $j$ be an admissible weight. Then both $L(\ell,j)$
and $L(\ell,j)^{c}$ are
irreducible weak $L(\ell,0)$-modules.
\end{coro}

\begin{rem} If $\ell$ is not a nonnegative integer, there are also other
types of irreducible
weak $L(\ell,0)$-modules. For a positive integral level $\ell$, it was
proved [DLiM2] that any weak module is completely reducible and any
irreducible weak $L(\ell,0)$-module is an irreducible integrable
highest weight $\hat\frak g$-module of level $\ell$.  This
distinguishes $L(\ell,0)$ for a positive integral level $\ell$ from
all the rational levels.
\end{rem}

\begin{rem}\label{r2.16}
It follows immediately from Propositions \ref{p2.8} and a complete reducibility
theorem
of Kac-Wakimoto (Theorem 4.1 of [KW2]) that any weak $L(\ell,0)$-module $M$
which is an $\hat\frak g$-module of level $\ell$ from the
category $\cal{O}$ is a direct sum of
irreducible modules $L(\ell,j)$ with admissible weight $j$.
\end{rem}

Nest, we shall prove a completely reducibility theorem for a category much
bigger than
the category $\cal{O}$. Recall the following theorem from [KK]:

\bt{tkk}
Let $\lambda, \mu\in H^{*}$. Then $L(\mu)$ is isomorphic to a subquotient
module of $M(\lambda)$ iff
the ordered pair $\{\lambda, \mu\}$ satisfies the following condition:
There exists a sequence $\beta_{1},\cdots,\beta_{k}$ of positive roots and a
sequence $n_{1},
\cdots,n_{k}$ of positive integers such that

(i) $\lambda-\sum_{i=1}^{k}n_{i}\beta_{i}=\mu$;

(ii)
$2(\lambda+\rho-n_{1}\beta_{1}-\cdots-n_{j-1}\beta_{j-1},\beta_{j})=n_{j}(\beta_{j},\beta_{j})$
for $1\le j\le k$.
\et

\begin{lem}\label{l4.12}
Let $\lambda, \mu$ be two distinct admissible weights. Then $L(\mu)$ is not
isomorphic to any
subquotient module of $M(\lambda)$.
\end{lem}

{\bf Proof.} Otherwise, by Theorem \ref{tkk}  we have
a sequence $\beta_{1},\cdots,\beta_{k}$ of positive roots and a sequence
$n_{1},
\cdots,n_{k}$ of positive integers satisfying (i)-(ii).
{}From [DGK] each $\beta_{i}$ is real.
 Then we obtain
\begin{eqnarray}
\<\mu+\rho,\beta_{k}\>=\<\lambda+\rho -\sum_{i=1}^{k}n_{i}\beta_{i},\beta_{k}\>
=\frac{2(\lambda+\rho-\sum_{i=0}^{k}n_{i}\beta_{i},\beta_{k})}{(\beta_{k},\beta_{k})}
=-n_{k}.
\end{eqnarray}
This contradicts the admissibility of $\mu$. $\;\;\;\;\Box$

\begin{lem}\label{l4.13}
Let $M$ be a weak $L(\ell,0)$-module such that $M$ is a highest weight
$\tilde{\frak{g}}$-module.
Then $M$ is irreducible.
\end{lem}

{\bf Proof.} Let $\lambda$ be the highest weight of $M$.
If $M$ contains a proper submodule $W$, there is a highest weight vector $u$ in
$W$ of weight
$\mu$ such that $\mu <\lambda$.Then both $\l$ and $\mu$ are admissible by
Corollary \ref{c2.8}. This contradicts Lemma \ref{l4.12}.
Then $M$ is irreducible.$\;\;\;\;\Box$

Recall from [K] that $\hat{\omega}$ is the involutory antiautomorphism  of
$\hat{\frak{g}}$,
which is
the negative Chevalley involution.
Let $M$ be a $\hat{\frak{g}}$-module of level $\ell$ such that $H$ local
finitely acts on $M$
with finite-dimensional
generalized $H$-eigenspaces. We define [DGK]
$M^{\hat{\omega}}=\oplus_{\lambda\in H^{*}}M_{\lambda}^{*}$ with the following
action
$(af)(u)=f(\hat{\omega}(a)u)$ for any $f\in M^{\hat{\omega}}, a\in
\hat{\frak{g}}, u\in M$. Then
$$(M^{\hat{\omega}})^{\hat{\omega}}\simeq M,\; L(\lambda)^{\hat{\omega}}\simeq
L(\lambda)
\;\;\;\mbox{ for any }\lambda\in H^{*}.$$

\bp{p4.14}
Let $\lambda_{1}, \lambda_{2}$ be admissible weights of level $\ell$. Then any
short exact
sequence
\begin{eqnarray}\label{es1}
0\rightarrow L(\lambda_{1})\rightarrow M\rightarrow L(\lambda_{2})\rightarrow 0
\end{eqnarray}
of weak $L(\ell,0)$-modules splits.
\ep

{\bf Proof.} First, since $H$ semisimply acts on $L(\lambda_{1})$ and
$L(\lambda_{2})$, $H$ acts
local finitely on $M$. Let $M=\oplus_{\lambda\in H^{*}}M_{\lambda}$
be the generalized $H$-eigenspace decomposition. Then the sequence (\ref{es1})
splits if and
only if
the following sequence splits:
\begin{eqnarray}\label{es2}
0\rightarrow L(\lambda_{2})\rightarrow M^{\hat{\omega}}\rightarrow
L(\lambda_{1})\rightarrow 0.
\end{eqnarray}
Without losing generality we may assume that $\lambda_{1}\not > \lambda_{2}$.
Let $u\in M_{\lambda_{2}}$ such that $u\notin L(\lambda_{1})$. Then
$N_{+}u\subseteq L(\lambda_{1})$. If $N_{+}u\ne 0$, there is a
$\beta\in {\Z}_{+}\alpha_{0}+{\Z}_{+}\alpha_{1}$ such that
$\lambda_{2}+\beta=\lambda_{1}$. This contradicts the assumption
$\lambda_{1}\not >\lambda_{2}$.
Thus $N_{+}u=0$. Set $U=U(\frak{g})u$. Let $W$ be the submodule generated by
$U$.
Since $L(\ell,U)$ as a $\hat{\frak{g}}$-module is isomorphic to some quotient
module of $W$,
$L(\ell,U)$ is a weak $L(\ell,0)$-module. From Proposition \ref{p2.8}, $H$
semisimplely acts on
$U$.
Then $u$ is a highest weight vector. By Lemma \ref{l4.13}, $W$ is irreducible.
Then we obtain
$M=W\oplus L(\lambda_{1})$. That is, sequence (\ref{es1}) splits.
$\;\;\;\;\Box$.

{}From Proposition \ref{p4.14} we have

\begin{coro}\label{c4.15}
Let $\lambda, \lambda_{1},\cdots, \lambda_{k}$ be admissible weights of level
$\ell$. Then any
short exact sequence
$$0\rightarrow L(\lambda_{1})\oplus \cdots \oplus L(\lambda_{k})\rightarrow
M\rightarrow L(\lambda)
\rightarrow 0$$
of weak $L(\ell,0)$-modules splits.
\end{coro}

\bt{t4.15}
Let $M$ be any weak $L(\ell,0)$-module
such that for any $u\in M$, there exists a positive integer $k$ such
that $(N_{+})^{k}u=0$. Then $M$ is a direct sum of
irreducible modules $L(\ell,j)$ with admissible weight $j$.
\et

{\bf Proof.} Set $\Omega(M)=\{m\in M|\frak g\otimes t\C[t]\cdot m=0\}.$ Then
the proof of Theorem 3.7 of [DLiM2] shows that $\Omega(M)\ne 0.$ Since
$e$ is locally nilpotent on $\Omega(M)$ we conclude that there exist vectors
$m\in M$ such that $N_+m=0.$  From the proof of Proposition \ref{p4.14}
we see that $H$ acts semisimply on $U(\frak g)m.$ Thus $M$ contains a
highest weight vector. It follows from Lemma \ref{l4.13} that
$M$ contains an
irreducible weak $L(\ell,0)$-module $L(\lambda)$. Let $W$ be the sum of all
irreducible weak
$L(\ell,0)$-submodules of $M$. We have to prove $M=W$. If $M\ne W$, there is a
submodule $E$ of $M$
such that $W\subseteq E,\; E/W\simeq L(\lambda)$ for some admissible weight
$\lambda$.
Let $u+W$ be a highest weight vector of $E/W$. Since $\hat{\frak{g}}$ is
finitely generated,
$$N_{+}u\subseteq L(\lambda_{1})\oplus L(\lambda_{2})\oplus \cdots \oplus
L(\lambda_{r})$$
for some $\lambda_{1},\cdots, \lambda_{r}$.
Set $W^{o}=L(\lambda_{1})\oplus L(\lambda_{2})\oplus \cdots \oplus
L(\lambda_{r})$. It follows from
Corollary \ref{c4.15} that the submodule generated by $u$ and $W^{o}$ is
completely reducible. Then
$u\in W$. This contradicts the assumption of $u$. Thus $M=W$. This finishes the
proof.
$\;\;\;\;\Box$

\begin{rem}\label{r2.17}
{}From Corollary \ref{c2.8} and Proposition \ref{t4.15} the set of equivalence
classes of
irreducible $L(\ell,0)$-modules
consists of $L(\ell,j)$ with $j\in {\Z}, 0\le j\le p-2$ and any (ordinary)
module is completely
reducible.
\end{rem}

\begin{rem}\label{r2.18}
In [Z], an associative algebra $A(V)$ was introduced for any vertex operator
algebra $V$ such that there is a natural one-to-one
correspondence between the set of equivalence classes of irreducible
${\Z}_{+}$-graded
weak $V$-modules and the set of equivalence classes of irreducible
$A(V)$-modules.
If $\ell$ is not a nonnegative integer, then $A(L(\ell,0))$ has
infinite-dimensional
irreducible modules so that $A(L(\ell,0))$ is infinite-dimensional. Therefore
(from [DLiM1]),
$L(\ell,0)$ is not rational. Because
Zhu's $C_{2}$-finiteness condition implies that
$A(L(\ell,0))$ is finite-dimensional, $L(\ell,0)$ does not satisfy the
$C_{2}$-finiteness condition.
\end{rem}

\section{${\Q}$-graded vertex operator algebras and the rationality of
$(L(\ell,0),\omega_{z})$}
If $j$ is not a nonnegative integer, homogeneous spaces of $L(\ell,j)$
are infinite-dimensional so that the character $\tr_{L(\ell,j)}
q^{L(0)}$ is not well-defined, where
$c_{\ell}=\frac{3\ell}{\ell+2}$. In [KW1]-[KW2], the modified characters
$\tr_{L(\ell,j)} q^{L(0)-{1\over 2}zh(0)-{c_{\ell}\over 24}}$ were
considered, where $z$ is a positive rational number less than $1.$
Noticing that
$L(0)-{1\over 2}zh(0)$ could be considered as the degree-zero
component of a Virasoro vector $\omega_{z}=\omega+{1\over 2}zh(-2){\bf
1}$ whose central charge is $c_{\ell,z}=c_{\ell}-6\ell z^{2}$, we
study $L(\ell,0)$ with respect to the new Virasoro element
$\omega_{z}$ in this section. We denote the new vertex operator algebra
by $\omega_{z}$ $(L(\ell,0),\omega_{z}).$
Note that $\omega_{z}$ $(L(\ell,0),\omega_{z})$ is ${\Q}$-graded instead of
${\Z}$-graded. This leads us to the study of ${\Q}$-graded vertex
operator algebras. In particular We extend Zhu's $A(V)$-theory and
Frenkel-Zhu's fusion rule formula to any ${\Q}$-graded vertex operator
algebra. That is, we construct an associative algebra $A(V)$ for any
${\Q}$-graded VOA $V$ and establish the one-to-one correspondence
between the set of equivalence classes of irreducible
${\Q}_{+}$-graded weak $V$-modules and the set of equivalence classes
of irreducible $A(V)$-modules.  If $V$ is ${1\over 2}\Z$-graded, our
construction $A(V)$ and related results coincide with those for vertex
operator superalgebra as developed in [KWa].
We also use complete
reducibility Theorem \ref{t4.15} to show that $(L(\ell,0),\omega_{z})$ is
rational.

A ${\Q}$-graded vertex operator algebra $V$ satisfies all the axioms
for a vertex operator algebra $V$ except that $V$ is ${\Q}$-graded by
weights instead of ${\Z}$-graded. In particular, a ${\Q}$-graded
vertex operator algebra is a generalized vertex operator algebra in
the sense of [DL].  The definitions of weak module and ordinary module
are as before. In the
definition of $\Z_+$-graded module for a $\Z$-graded vertex operator
algebra, replacing $\Z$ by $\Q$ gives a $\Q_+$-graded module for a
$\Q$-graded vertex operator algebra.

\begin{de}\label{d3.4}
A ${\Q}$-graded vertex operator algebra $V$ is called {\em rational} if any
${\Q}_{+}$-graded
weak $V$-module is completely reducible.
\end{de}

Let $V=\oplus_{\alpha\in
{\Q}}V_{\alpha}$ be a ${\Q}$-graded vertex operator algebra.
Then $V_{{\Z}}=\oplus_{n\in {\Z}}V_{n}$ is a ${\Z}$-graded
(ordinary) vertex operator algebra.
Just as in [FFR] and [Li2], one obtains a
${\Q}$-graded Lie algebra $G(V)=\oplus_{\alpha\in {\Q}}G(V)_{\alpha}$ as the
quotient space of ${\C}[x,x^{-1}]\otimes V$ modulo
$({d\over dx}\otimes 1+1\otimes L(-1))({\C}[x,x^{-1}]\otimes V).$ Here the
Lie bracket is induced from
$$[x^m\otimes u, x^n\otimes]=
\sum_{i=0}^{\infty}{m\choose i}x^{m+n-i}\otimes u_iv$$
for $u,v\in V$ and the degree of $x^n\otimes u+({d\over dx}\otimes 1+1\otimes
L(-1))({\C}[x,x^{-1}]\otimes V)$ is $\wt u-n-1$ for homogeneous $u.$
Set \begin{eqnarray}
G(V)_{\pm}=\oplus_{\alpha >0}G(V)_{\pm \alpha}.
\end{eqnarray}

Let $U$ be any $G(V)_{0}$-module.
Then we form the following induced module:
\begin{eqnarray}
M(U)=U(G(V))\otimes_{U(G(V)_0\oplus G(V)_{-})}U
\end{eqnarray}
where $G(V)_-$ acts trivially on $U.$
Then $M(U)$ is a lower-truncated ${\Q}$-graded $G(V)$-module
generated by the lowest-degree subspace $U$. Let $U^*$ be
the dual space of $U$ and extend $U^*$ to $M(U)$ by letting $U^*$ annihilate
$\oplus_{n>0}M(U)(n).$ We denote such a pair by
$\langle u', v\rangle$
for $u'\in U^{*}$ and $v\in M(U).$  Set
\begin{eqnarray}
I=\{v\in M(U)|\langle u', av\rangle =0\;\;\mbox{for any }u'\in U^{*},
a\in U(G(V))\}.
\end{eqnarray}
Then it is clear that $I$ is a $G(V)$-submodule of $M(U)$. Let $L(U)$
be the quotient module of $M(U)$ modulo $I$.

Let $V$ be a ${\Q}$-graded vertex operator algebra. First we define a function
$\varepsilon$ for all homogeneous elements of $V$ as follows:
$\varepsilon (a)=1$ if ${\rm wt}a\in {\Z}$, $\varepsilon (a)=0$ if ${\rm
wt}a\notin {\Z}$.
For any homogeneous element $a\in V$, we define:
\begin{eqnarray}\label{d3.5}
a*b=\varepsilon (a){\rm Res}_{x}\frac{(1+x)^{[{\rm wt}a]}}{x}Y(a,x)b\;\;\mbox{
 for any }b\in V,
\end{eqnarray}
where $[\cdot ]$ denotes the greatest-integer function.
Then extend ``$*$'' to a bilinear product on $V$. Let $O(V)$ be the subspace
of $V$ linearly spanned by
\begin{eqnarray}
{\rm Res}_{x}\frac{(1+x)^{[{\rm wt}a]}}{x^{1+\varepsilon (a)}}Y(a,x)b
\end{eqnarray}
for any homogeneous element $a\in V$ and for any $b\in V$.
Using $(1+x)^{m}=\sum_{i=0}^{m}{m\choose i}x^{i}$ one can prove
\begin{eqnarray}
{\rm Res}_{x}\frac{(1+x)^{[\alpha]+m}}{x^{n+1+\varepsilon (a)}}Y(a,x)b\in I
\end{eqnarray}
for $n\ge m\ge 0$. Let $M$ be any weak
$V$-module. Then we define
$$\Omega(M)=\{ u\in M| a_{m}u=0\;\;\;\mbox{for }a\in V, m>{\rm wt}a-1\}.$$
Define $o$ to be the linear map from $V$ to ${\rm End}\Omega (M)$ such that
$o(a)=\varepsilon (a)a_{[{\rm wt}a]-1}$ for any homogeneous element $a$ of
$V$. Generalizing Theorems 2.1.1 and 2.1.2 of Zhu we obtain

\bt{tfz} (a) The subspace $O(V)$ is a two-sided ideal of $V$ with respect to
the product ``$*$'' and
$A(V)=V/O(V)$ is an associative algebra with identity ${\bf 1} +O(V).$
Moreover, $\o+O(V)$ lies in the center of $A(V).$

(b) For any weak $V$-module $M$, $\Omega(M)$ is an $A(V)$-module with $a$
acts as $o(a)$.
\et

The proof is the same as in the twisted case (see the proofs of Proposition
2.3 and Theorem 5.3 in [DLiM1]).

Similarly, for a weak $V$-module $M$, we define $O(M)$ to be the subspace of
$M$
linearly spanned by
\begin{eqnarray}
{\rm Res}_{x}\frac{(1+x)^{[{\rm wt}a]}}{x^{1+\varepsilon (a)}}Y(a,x)u
\end{eqnarray}
for any homogeneous element $a\in V$ and for any $u\in M$. The following
theorem is an analogue of Theorems 1.5.1 and 1.5.2
 of [FZ] (also see [KWa] and [Li2]):
\bt{tz} (a) The quotient space $A(M)= M/O(M)$ is an $A(V)$-bimodule with
the following left and right actions:
\begin{eqnarray}
& &a*u=\varepsilon (a){\rm Res}_{x}\frac{(1+x)^{[{\rm wt}a]}}{x}Y(a,x)u,\\
& &u*a=\varepsilon (a){\rm Res}_{x}\frac{(1+x)^{[{\rm wt}a]-1}}{x}Y(a,x)u
\end{eqnarray}
for any homogeneous $a\in V$ and for any $ u\in M$.

(b) Let $W_{1}, W_{2}, W_{3}$ be irreducible $V$-modules and suppose $V$ is
rational\footnote{It
was pointed out in [Li2] that this condition is necessary and a proof was
supplied}.
Then there is a linear isomorphism from the space
${\rm Hom}_{A(V)}(A(W_{1})\otimes_{A(V)}W_{2}(0), W_{3}(0))$
to the space of intertwining operators of type
$\left(\begin{array}{c}W_{3}\\W_{1} W_{2}\end{array}\right)$.
\et

{\bf Proof.} Let $I(M)$ be the subspace of $O(M)$ linearly spanned by
$${\rm Res}_{x}\frac{(1+x)^{{\rm wt}a}}{x^{2}}Y(a,x)u$$
for any homogeneous element $a\in V_{{\Z}}$ and for any $u\in M$.
Then $A_{V_{{\Z}}}(M)=M/I(M)$ is the $A(V_{{\Z}})$-bimodule defined in [FZ].
Thus it is enough for us to prove that the subspace $O(M)/I(M)$ is
a  sub-bimodule of $A_{V_{\Z}}(M)$.
Since the proof this is parallel to the proof of
Theorem \ref{tfz} we omit the proof.

The proof of (b) is similar to that for the $\Z$-graded vertex operator
algebra as in [Li2].$\;\;\;\;\Box$

By definition  $A(V)$ is a quotient algebra of $A(V_{{\Z}})$.
It is clear that $A(V)_{Lie}$ is a quotient algebra of
$G(V)_{0}$. Then for any $A(V)$-module $U$, we may naturally view $U$
as a $G(V)_{0}$-module.

\begin{prop}\label{p3.3} For any $A(V)$-module $U$, $L(U)$ is a
weak $V$-module.
\end{prop}

{\bf Proof.} The proof is the same as in the ordinary case
(see [Li2] and [Z]) or the twisted case (see the proof of Theorem
6.3 of [DLiM1]).$\;\;\;\;\Box$

The following is a  gneralization of Theorem 2.2.2 of [Z].
See [Li2] or [DLiM1] for a similar proof.

\bt{t3.4}
The functor $\Omega$ gives rise to a one-to-one correspondence between the set
of equivalence
classes of irreducible ${\Q}_{+}$-graded weak $V$-modules and the set of
equivalence classes of irreducible $A(V)$-modules.
\et

As in the case of $\Z$-graded vertex operator algebra, we have
(see the proof of Theorem 8.1 of [DLiM1]):

\begin{prop}\label{p3.5}
If $V$ is rational, $A(V)$ is semisimple and any ${\Q}_{+}$-graded weak
$V$-module is a
direct sum of irreducible ordinary $V$-modules.
\end{prop}

Let $V$ be a ${\Q}$-graded vertex operator algebra. Then it is clear that
$\exp (2\pi i L(0))$ is an automorphism of $V$.
Let $M=\oplus_{h\in
{\C}}M_{(h)}$ be a $V$-module. Following [FHL], let
$M'=\oplus_{h\in {\C}}M_{h}^{*}$ be the restricted dual of $M$ and define
\begin{eqnarray}
\langle Y(a,x)f,u\rangle=\langle f, Y(e^{xL(1)}(e^{\pi
i}x^{-2})^{L(0)}a,x^{-1})u
\rangle
\end{eqnarray}
for any $f\in M', a\in V, u\in M$.
The following proposition is essentially  proved in [Li3].

\begin{prop}\label{p3.6} The pair $(M', Y(\cdot,x))$ gives rise to a
$\sigma^{2}$-twisted
$V$-module, where $\sigma =\exp (2\pi i L(0)).$
\end{prop}

\begin{rem}\label{r3.7} If $V$ is ${1\over 2}{\Z}$-graded, then $M'$ is a
$V$-module because $\sigma^{2}= id_{V}$. Therefore, we obtain a new functor
from $V$-modules to $V$-modules. It is important to notice that the vertex
operator algebra $V$ may not be isomorphic to its own contragredient dual.
\end{rem}

Let $V$ be a ${\Q}$-graded vertex operator algebra and let $M$ be any weak
$V$-module. Define
$C_{2}(M)$ to be the subspace
linearly spanned by $a_{-2}M$ for $a\in V_{{\Z}}$ and by $a_{-1}M$ for
$a\in V_{n}, n\not \in {\Z}$.
Define bilinear products ``$\cdot$'' and ``$\circ$'' on $V$ as follows:
For $a\in V_{m},b\in V_{n}$ we define
$a\cdot b=a_{-1}b$ and $a\circ b=a_{0}b$ if $m,n\in {\Z}$, otherwise we define
$a\cdot b=0$
and $a\circ b=0$.

\begin{lem}\label{l3.9}
The defined subspace $C_{2}(V)$ is a two-sided ideal for both $(V, \cdot)$ and
$(V, \circ)$.
\end{lem}

{\bf Proof.} Let $a\in V_{m}, b\in V_{n}, c\in V_{k}$. If $m\not \in {\Z}$ or
$n+k\not \in {\Z}$, by definition we have:
$a\cdot b_{-r}c=0$ and $a\circ b_{-r}c=0$ for $r=1$ or $2$. If $m, n+k\in
{\Z}$, we get
$$a_{-j}(b_{-r}c)=b_{-r}a_{-j}c
+\sum_{i=0}^{\infty}\left(\begin{array}{c}-j\\i\end{array}\right)(a_{i}b)_{-j-r-i}c$$
for $j=-1$ or $0$.
Then $a\cdot (b_{-r}c), a\circ (b_{-r}c)\in C_{2}(V)$. Then the proof is
complete.
$\;\;\;\;\Box$

Set $A_{2}(M)=M/C_{2}(M)$.
By Lemma \ref{l3.9} we obtain a quotient algebra $A_{2}(V)=V/C_{2}(V)$.
Similarly
to [Z] we have:

\begin{prop}\label{p3.10}
The quotient algebra $(A_{2}(V),\cdot)$ is a
commutative associative algebra with the vacuum vector
${\bf 1}$ as its identity and $(A_{2}(V),\circ)$ is a Lie algebra such that
$$(a\cdot b)\circ c=a\cdot (b\circ c)+(a\circ c)\cdot b$$
for any $a,b,c\in A_{2}(V)$. Therefore $(A_{2}(V), \cdot, \circ)$ is a Poisson
Lie algebra.
\end{prop}

\begin{de}\label{d3.11}
If $A_{2}(V)$ is finite-dimensional, we say $V$ is {\em $C_{2}$-finite} or $V$
satisfies the $C_{2}$-finiteness condition. If $V$ as a Virasoro algebra module
is generated
by primary vectors, we say that $V$ satisfies the {\em Virasoro condition}. If
$V$ as a vertex
operator algebra is generated by $\omega$ and all primary vectors, we say that
$V$ satisfies
the {\em primary-field condition}.
\end{de}

\begin{rem}\label{r3.13} It has been proved in [Z] that
if $V$ is a rational vertex operator algebra with integral weights satisfying
the
$C_{2}$-finiteness condition and the Virasoro condition,
then the space linearly spanned by $\tr_{M}q^{L(0)-{c\over 24}}$,
where $M$ runs through all irreducible $V$-modules, is modular invariance.
If one replaces the Virasoro condition by the primary-field condition,
one can check that Zhu's theorem also holds.
\end{rem}

Recall the following proposition from [DLinM].

\begin{prop}\label{DLiM} Let $(V,Y,{\bf 1},\omega)$ be a vertex
operator algebra of rank $r$ and let
$h\in V$ satisfying the following conditions:
\begin{eqnarray}
L(n)h=\delta_{n,0}h,\;
h_{n}h=\delta_{n,1}\lambda {\bf 1}\;\;\mbox{{\it for }}n\in {\Z}_{+},
\end{eqnarray}
where $\lambda$ is a complex number. Then $(V,Y,{\bf 1},\omega+h_{-2}{\bf 1})$
is a vertex
algebra of rank $r -12\lambda$.
\end{prop}

Now we go back to the vertex operator algebra $L(\ell,0)$. For any $z\in {\Q}$,
we
set $\omega_{z}=\omega+{1\over 2}zh(-2){\bf 1}$.
Then it follows from Proposition 4.1 of
[DLinM] that $\omega_{z}$ is a new Virasoro vector of
$L(\ell,0)$ with a central charge $\frac{3\ell}{\ell+2}-6\ell z^{2}$.
Thus $L_{z}(0)=(\omega_{z})_{1}=L(0)-{1\over 2}zh(0)$ so that
\begin{eqnarray}
& &[L_{z}(0), x^{m}\otimes h]=-m(x^{m}\otimes h);\\
& &[L_{z}(0), x^{m}\otimes e]=(-m-z)(x^{m}\otimes e);\\
& &[L_{z}(0), x^{m}\otimes f]=(-m+z)(x^{m}\otimes f)
\end{eqnarray}
for any $m\in {\Z}$. In general, $V$ is ${\Q}$-graded by weights
with respect to $L_{z}(0)=L(0)-{1\over 2}zh(0)$ instead of ${\Z}$-graded.
Consequently, we obtain a ${\Q}$-grading for
$\tilde{sl}_{2}$ satisfying the conditions:
\begin{eqnarray}
\deg (x^{n}\otimes e)=-n-z;\;\deg (x^{n}\otimes f)=-n+z,\;\deg (x^{n}\otimes
h)=-n \;\;\;
\mbox{ for }n\in {\Z}.
\end{eqnarray}
For a positive integral level $\ell$, all irreducible $L(\ell,0)$-modules are
integral
modules. For a general rational level $\ell$, admissible weight $j$ may be
non-integral.
To make the graded spaces of $L(\ell,j)$ be finite-dimensional, we assume $z\in
{\Q}, 0<z<1$.

Let $M=\oplus_{n\in {\Q}_{+}}M(n)$ be any ${\Q}_{+}$-graded weak
$(L(\ell,0),\omega_{z})$-module.
Since $x^{n}\otimes e, x^{n+1}\otimes f, x^{n+1}\otimes h$ for $n\in {\Z}_{+}$
have negative
degrees
with respect to the operator $L_{z}(0)$, it is clear that $M$ satisfies the
condition of
Proposition \ref{t4.15} so that $M$ is completely reducible. Then we obtain

\bt{t4.6}
The ${\Q}$-graded vertex operator algebra $(L(\ell,0),\omega_{z})$ is rational
and all irreducible modules (up to equivalence) are $L(\ell,j)$ for the
admissible weights $j.$
\et

\begin{rem}\label{r3.2} It is easy to check that each eigenspace for
$L_z(0)$ in $L(\ell,j)$ is finite-dimensional for admissible weight $j.$
Thus $\tr_{L(\ell,j)}q^{L_z(0)}=\tr_{L(\ell,j)}q^{L(0)-{1\over 2}zh(0)}$
is well defined and is equal to $\tr_{L(\ell,j)^{c}}q^{L_z(0)}.$
In fact they are convergent in upper half plane (see Section 5).
\end{rem}

\begin{rem}\label{r3.12}
Since $L_{z}(n)\!=\!L(n)-{1\over 2}(n+1)zh(n)$, $e,f$ are primary vectors in
$(\!L(\ell,0),\omega_{z}\!)$.
Because $L(\ell,0)$ as a vertex operator algebra is generated by $e,f$,
$(L(\ell,0),\omega_{z})$
satisfies the primary-field condition.
\end{rem}

\section{Fusion rules and $C_{2}$-finiteness of  $(L(\ell,0), \omega_{z})$}
The main goal of this section is to calculate the fusion rules and prove that
$(L(\ell,0), \omega_{z})$ is
$C_{2}$-finite. Throughout this section we assume that $\ell=-2+{p\over q}$,
where $p$ and
$q$ are coprime positive integers with $p\ge 2$ and that $z$ is a fixed
rational number
satisfying
$0<z<1$ (under which  certain traces converge in some domain [KW1-2]).

Let $M$ be any weak $V(\ell,\C)$-module. Since
\begin{eqnarray}\label{d4.1}
{\rm wt}h(-1)=1,\;{\rm wt }e(-1)=1-z,\;{\rm wt}f(-1)=1+z,
\end{eqnarray}
we have
\begin{eqnarray}
{\rm Res}_{x}\frac{(1+x)^{[{\rm wt}f]}}{x^{m}}Y(f,x)u
&=&(f(-m)+f(1-m))u;\label{d4.2}\\
{\rm Res}_{x}\frac{(1+x)^{[{\rm wt}e]}}{x^{m}}Y(e,x)u
&=&e(-m)u;\\
{\rm Res}_{x}\frac{(1+x)^{{\rm wt}h}}{x^{m+1}}Y(h,x)u
&=&(h(-m-1)+h(-m))u\label{d4.4}
\end{eqnarray}
for any positive integer $m$ and for $u\in M$. By definition all those
elements in (\ref{d4.2})-(\ref{d4.4}) are in $O(M)$.

\begin{prop}\label{p4.1} Let $M$ be any weak $V(\ell,\C)$-module. Then
the space $O(M)$ is spanned by the
all the elements in  (\ref{d4.2})-(\ref{d4.4}).
\end{prop}

{\bf Proof.} Let $W$ be the subspace linearly spanned by the all
elements in  (\ref{d4.2})-(\ref{d4.4}). Set
\begin{equation}\label{dd4.4}
C={\C}[x^{-1}](x^{-1}+1)\otimes f+{\C}[x^{-1}]x^{-1}\otimes e
+{\C}[x^{-1}](x^{-2}+x^{-1})\otimes h.
\end{equation}
Then $W=C\cdot M.$
Since $[h(-k),C]\subseteq C$ for any positive integer $k$,
we get $h(-k)W\subseteq W$.

Let $L$ be the linear span of homogeneous
elements $a$ of $V(\ell,\C)$ such that for any positive integer $n$,
\begin{eqnarray}
{\rm Res}_{x}\frac{(1+x)^{[{\rm wt}a]}}{x^{n+\varepsilon(a)}}Y(a,x)M
\subseteq W.
\end{eqnarray}
We shall prove that  $L$ is equal to $V(\ell,\C).$

For any homogeneous element $a$ of $L$ and for any nonnegative
integers $m\ge n$, we have
\begin{eqnarray}\label{e4.6}
{\rm Res}_{x}\frac{(1+x)^{[{\rm wt}a]+n}}{x^{m+1+\varepsilon (a)}}Y(a,x)M
\subseteq W
\end{eqnarray}
because $\displaystyle{\frac{(1+x)^{[{\rm wt}a]+n}}{x^{m+1+\varepsilon(a)}}=
\sum_{i=0}^{\infty}{n\choose i}
\frac{(1+x)^{[{\rm wt}a]}}{x^{m-i+1+\varepsilon(a)}}}$.

Let $a$ be any homogeneous element of $L$ and let
$k$ be any positive integer. Then for any $n\in {\N},u\in M$, we have:
\begin{eqnarray}
& &{\rm Res}_{z_{2}}\frac{(1+z_{2})^{[{\rm wt}h(-k)a]}}{z_{2}^{n+\varepsilon
(h(-k)a)}}
Y(h(-k)a,z_{2})u\nonumber\\
&=&{\rm Res}_{z_{0}}{\rm Res}_{z_{2}}\frac{(1+z_{2})^{[{\rm wt}a]+k}}
{z_{2}^{n+\varepsilon(a)}}z_{0}^{-k}Y(Y(h,z_{0})a,z_{2})u\nonumber\\
&=&{\rm Res}_{z_{1}}{\rm Res}_{z_{2}}
\frac{(1+z_{2})^{[{\rm
wt}a]+k}}{z_{2}^{n+\varepsilon(a)}}(z_{1}-z_{2})^{-k}Y(h,z_{1})
Y(a,z_{2})u\nonumber\\
& &-{\rm Res}_{z_{2}}{\rm Res}_{z_{1}}
\frac{(1+z_{2})^{[{\rm
wt}a]+k}}{z_{2}^{n+\varepsilon(a)}}(-z_{2}+z_{1})^{-k}Y(a,z_{2})
Y(h,z_{1})u\nonumber\\
&=&{\rm Res}_{z_{2}}\sum_{i=0}^{\infty}{-k\choose i}
(-z_{2})^{i}\frac{(1+z_{2})^{[{\rm wt}a]+k}}{z_{2}^{n+\varepsilon (a)}}
h(-k-i)Y(a,z_{2})u\nonumber\\
& &-{\rm Res}_{z_{2}}\sum_{i=0}^{\infty}{-k\choose i}
(-1)^{k+i}\frac{(1+z_{2})^{[{\rm wt}a]+k}}{z_{2}^{n+k+i+\varepsilon(a)}}
Y(a,z_{2})h(i)u\nonumber\\
&\equiv&{\rm Res}_{z_{2}}\sum_{i=0}^{\infty}{-k\choose i}
z_{2}^{i}\frac{(1+z_{2})^{[{\rm wt}a]+k}}{z_{2}^{n+\varepsilon (a)}}
h(-k)Y(a,z_{2})u \;\;\mbox{ mod }W\nonumber\\
&=&h(-k){\rm Res}_{z_{2}}\frac{(1+z_{2})^{[{\rm
wt}a]}}{z_{2}^{n+\varepsilon(a)}}
Y(a,z_{2})u \nonumber\\
&\equiv&0\;\;\;\mbox{ mod }W.
\end{eqnarray}
Here we used the relation $h(-k-i)w\equiv (-1)^ih(-k)w$ (mod $M$)  which
follows from (\ref{d4.4}).
Therefore $h(-k)L\subseteq L$ for any $k\in {\N}$.

Similarly, we have
\begin{eqnarray}
& &{\rm Res}_{z_{2}}\frac{(1+z_{2})^{[{\rm wt}e(-k)a]}}{z_{2}^{n+\varepsilon
(e(-k)a)}}
Y(e(-k)a,z_{2})u\nonumber\\
&=&{\rm Res}_{z_{2}}\sum_{i=0}^{\infty}{-k\choose i}
(-z_{2})^{i}\frac{(1+z_{2})^{[{\rm wt}e(-k)a]}}{z_{2}^{n+\varepsilon(e(-k)a)}}
e(-k-i)Y(a,z_{2})u\nonumber\\
& &-{\rm Res}_{z_{2}}\sum_{i=0}^{\infty}{-k\choose i}
(-1)^{k+i}\frac{(1+z_{2})^{[{\rm
wt}e(-k)a]}}{z_{2}^{n+k+i+\varepsilon(e(-k)a)}}
Y(a,z_{2})e(i)u\nonumber\\
&\equiv&-{\rm Res}_{z_{2}}\sum_{i=0}^{\infty}{-k\choose i}
(-1)^{k+i}\frac{(1+z_{2})^{[{\rm
wt}e(-k)a]}}{z_{2}^{n+k+i+\varepsilon(e(-k)a)}}
Y(a,z_{2})e(i)u \;\;\mbox{ mod }W.\label{d4.8}
\end{eqnarray}
If ${\rm wt}a\in {\Z}$, then $[{\rm wt}e(-k)a]={\rm wt}a+k-1$ and
$\varepsilon(e(-k)a)=0$.
Then the last formula in (\ref{d4.8}) is equal to
\begin{eqnarray*}
-{\rm Res}_{z_{2}}\sum_{i=0}^{\infty}{-k\choose i}
(-1)^{k+i}\frac{(1+z_{2})^{{\rm wt}a+k-1}}{z_{2}^{n+k+i}}
Y(a,z_{2})e(i)u
\end{eqnarray*}
which is in $W$ by (\ref{e4.6}) as ${\rm wt}e=1-z<1$ and
$[{\rm wt}a]+k-1\le [{\rm wt}e(-k)a]\le [{\rm wt}a]+k$.
A similar discussion using (\ref{e4.6}) shows that the last expression
of (\ref{d4.8}) is also in $W$ if ${\rm wt}a\notin {\Z}.$
Thus
${\rm Res}_{z_{2}}\frac{(1+z_{2})^{[{\rm wt}e(-k)a]}}{z_{2}^{n+\varepsilon
(e(-k)a)}}
Y(e(-k)a,z_{2})u\in W$.

Analogously,
\begin{eqnarray}
& &{\rm Res}_{z_{2}}\frac{(1+z_{2})^{[{\rm wt}f(-k)a]}}{z_{2}^{n+\varepsilon
(f(-k)a)}}
Y(f(-k)a,z_{2})u\nonumber\\
&=&{\rm Res}_{z_{2}}\sum_{i=0}^{\infty}{-k\choose i}
(-z_{2})^{i}\frac{(1+z_{2})^{[{\rm wt}a+z]+k}}{z_{2}^{n+\varepsilon (f(-k)a)}}
f(-k-i)Y(a,z_{2})u\nonumber\\
& &-{\rm Res}_{z_{2}}\sum_{i=0}^{\infty}{-k\choose i}
(-1)^{k+i}\frac{(1+z_{2})^{[{\rm wt}a+z]+k}}{z_{2}^{n+k+i+\varepsilon
(f(-k)a)}}
Y(a,z_{2})f(i)u\nonumber\\
&\equiv& {\rm Res}_{z_{2}}\sum_{i=0}^{\infty}{-k\choose i}(-1)^{k}
z_{2}^{i}\frac{(1+z_{2})^{[{\rm wt}a+z]+k}}{z_{2}^{n+\varepsilon (f(-k)a)}}
f(0)Y(a,z_{2})u\nonumber\\
& &-{\rm Res}_{z_{2}}\sum_{i=0}^{\infty}{-k\choose i}
(-1)^{k+i}\frac{(1+z_{2})^{[{\rm wt}a+z]+k}}{z_{2}^{n+k+i+\varepsilon
(f(-k)a)}}
Y(a,z_{2})f(i)u\;\;\mbox{ mod }W\nonumber\\
&\equiv& {\rm Res}_{z_{2}}(-1)^{k}
\frac{(1+z_{2})^{[{\rm wt}a+z]}}{z_{2}^{n+\varepsilon (f(-k)a)}}
f(0)Y(a,z_{2})u\nonumber\\
& &-{\rm Res}_{z_{2}}
(-1)^{k}\frac{(1+z_{2})^{[{\rm wt}a+z]+k}}{z_{2}^{n+k+\varepsilon (f(-k)a)}}
Y(a,z_{2})f(0)u\nonumber\\
& &-{\rm Res}_{z_{2}}\sum_{i=1}^{\infty}{-k\choose i}
(-1)^{k+i}\frac{(1+z_{2})^{[{\rm wt}a+z]+k}}{z_{2}^{n+k+i+\varepsilon
(f(-k)a)}}
Y(a,z_{2})f(i)u\;\;\mbox{ mod }W\nonumber\\
&\equiv& {\rm Res}_{z_{2}}(-1)^{k}
\frac{(1+z_{2})^{[{\rm wt}a+z]}}{z_{2}^{n+\varepsilon (f(-k)a)}}
(f(0)Y(a,z_{2})-Y(a,z_{2})f(0))u\nonumber\\
& &-{\rm Res}_{z_{2}}\sum_{i=1}^{k}{k\choose i}
(-1)^{k}\frac{(1+z_{2})^{[{\rm wt}a+z]}}{z_{2}^{n+i+\varepsilon (f(-k)a)}}
Y(a,z_{2})f(0)u\nonumber\\
& &-{\rm Res}_{z_{2}}\sum_{i=1}^{\infty}{-k\choose i}
(-1)^{k+i}\frac{(1+z_{2})^{[{\rm wt}a+z]+k}}{z_{2}^{n+k+i+\varepsilon
(f(-k)a)}}
Y(a,z_{2})f(i)u\nonumber\\
&\equiv& {\rm Res}_{z_{2}}(-1)^{k}
\frac{(1+z_{2})^{[{\rm wt}f(0)a]}}{z_{2}^{n+\varepsilon (f(0)a)}}
Y(f(0)a,z_{2})u\nonumber\\
& &-{\rm Res}_{z_{2}}\sum_{i=1}^{k}{k\choose i}
(-1)^{k}\frac{(1+z_{2})^{[{\rm wt}a+z]}}{z_{2}^{n+i+\varepsilon (f(-k)a)}}
Y(a,z_{2})f(0)u\nonumber\\
& &-{\rm Res}_{z_{2}}\sum_{i=1}^{\infty}{-k\choose i}
(-1)^{k+i}\frac{(1+z_{2})^{[{\rm wt}a+z]+k}}{z_{2}^{n+k+i+\varepsilon
(f(-k)a)}}
Y(a,z_{2})f(i)u\;\;\mbox{ mod }W.\label{d4.9}
\end{eqnarray}
Since $\deg f(0)a=\deg a$, by the induction hypothesis, we have
$${\rm Res}_{z_{2}}\frac{(1+z_{2})^{[{\rm wt}f(0)a]}}{z_{2}^{n+\varepsilon
(f(0)a)}}
Y(f(0)a,z_{2})u\in W.$$
Notice that $[{\rm wt}a+z]=[{\rm wt}a]$, or $[{\rm wt}a]+1$. If $[{\rm
wt}a+z]=[{\rm wt}a]$, by (\ref{e4.6}) the last two terms in (\ref{d4.9}) are in
$W$
no matter what $\varepsilon (f(-k)a)$ is.
 If $[{\rm wt}a+z]=[{\rm wt}a]+1$, then
${\rm wt}a\notin {\Z}$ so that $\varepsilon (a)=0$,
it is clear that the last two terms in (\ref{d4.9}) are in $W$ again by
(\ref{e4.6}).

Since ${\bf 1}\in L$ and $V(\ell,\C)=U(x^{-1}{\C}[x^{-1}]\otimes \frak g){\bf
1}$ we get $L=V(\ell,\C)$. Thus
$O(M)\subseteq W$. Therefore $O(M)=W.\;\;\;\;\Box$

\begin{prop}\label{p4.2} The associative algebra $A(V(\ell,\C))$ for
${\Q}$-graded vertex
operator algebra $(V(\ell,\C),\omega_{z})$ is isomorphic to the
polynomial algebra ${\C}[x]$.
\end{prop}

{\bf Proof.} Define a linear map $\psi$ from ${\C}[x]$ to $A(V(\ell,\C))$ as
follows
\begin{eqnarray}
\psi (g(x))=g(h(-1)){\bf 1}+O(V(\ell,\C))
\end{eqnarray}
for $g(x)\in {\C}[x]$. Since $[h(-1),h(0)]\!=\!0$ and $h(0){\bf 1}\!=\!0$, we
get
$g(h(-1)){\bf 1}=g(h(-1)+h(0)){\bf 1}$ for any $g(x)\in {\C}[x]$.
Since wt\,$h(-1)=1$ it follows from definition (\ref{d3.5}) that
$h(-1)*u=(h(-1)+h(0))u$ for $u\in V(\ell,\C).$
Thus $\psi$ is an algebra homomorphism.
Recall $N_-$ and $C$ from (\ref{d2.3}) and (\ref{dd4.4}). Then
$N_{-}=B\oplus {\C}f(0)\oplus {\C}h(-1).$ We have
$U(N_{-})=U(C)U({\C}h(-1))U({\C}f(0))$. By Proposition \ref{p4.1}
$$O(V(\ell,\C))=CV(\ell,\C)= CU(N_{-}){\bf 1}\simeq CU(C)U({\C}h(-1)).$$
Therefore $\psi$ is an isomorphism.$\;\;\;\;\Box$

\begin{prop}\label{p4.3} The $A(V)$-bimodule $A(M(\ell,j))$ is isomorphic
to ${\C}[x,y]$ with the bi-action as follows:
\begin{eqnarray}\label{d4.11}
x* f(x,y)=(x+j-2y{\partial\over \partial y})f(x,y)  ,\; f(x,y)* x=xf(x,y)
\end{eqnarray}
for any $f(x,y)\in {\C}[x,y]$.
\end{prop}

{\bf Proof.} Let $v$ be a (nonzero) lowest weight vector of $M(\ell,j)$.
Then as in the proof of Proposition \ref{p4.2} we have
$$O(M(\ell,j))=CU(C)U({\C}h(-1))U({\C}f(0))v\simeq
CU(C)U({\C}h(-1))U({\C}f(0)).$$
Then
$$A(M(\ell,j))=\oplus_{m,n\in {\Z}_{+}}{\C}(h(-1)^{m}f(0)^{n}+O(M(\ell,j)).$$
By the definition of the left and right actions of $A(V(\ell,\C))$
on $A(M(\ell,j))$ in Theorem \ref{tz}, we have
\begin{eqnarray}\label{e4.15}
h(-1)* (h(-1)^{m}f(0)^{n}v)&=&(h(-1)+h(0))h(-1)^{m}f(0)^{n}v\nonumber\\
&=&(h(-1)+j-2n)h(-1)^{m}f(0)^{n}v
\end{eqnarray}
and
\begin{eqnarray}
(h(-1)^{m}f(0)^{n}v)* h(-1)=h(-1)(h(-1)^{m}f(0)^{n}v)=h(-1)^{m+1}f(0)^{n}v.
\end{eqnarray}
The proposition follows immediately if we set $x=h(-1)+O(M(\ell,j)),
y=f(0)+O(M(\ell,j))$.
$\;\;\;\;\Box$

As a corollary of Propositions \ref{p2.8}, \ref{p3.5} and Theorem \ref{t3.4}
we obtain

\begin{coro}\label{c4.5}
The associative algebra $A(L(\ell,0))$ is semisimple and
isomorphic to the quotient algebra ${\C}[x]/\langle f(x)\rangle$
of the polynomial algebra ${\C}[x]$ in $x$, where
\begin{eqnarray}
f(x)=\prod_{r=0}^{p-2}\prod_{s=0}^{q-1}(x-r+st).
\end{eqnarray}
\end{coro}

The following lemma is useful for calculating $A(L(\ell,j)).$ The reader
can refer to [FZ] for a proof.

\begin{lem}\label{l4.4} (a) Let $V$ be a vertex operator algebra
and let $M$ be a $V$-module with a submodule $W$. Set $\bar{M}=M/W$.
Then as an $A(V)$-bimodule $A(\bar{M})\simeq M/(O(M)+W)$.

(b) If $I$ is an ideal of $V$ then $(I+O(V))/0(V)$ is a 2-sided ideal of $A(V)$
and $A(V/I)$ is isomorphic to
$A(V)/((I+O(V))/O(V)).$

(c)  If $I$ is an ideal of $V$,  and $I\cdot W \subset M$ ($ I \cdot W$ means
the linear span of elements $ v_nw$ for $v\in I, n\in{\bf Z}$ and $w\in
W$), then $I*A(M)  \subset (W+O(M))/O(M),$
$A(M)*I \subset (W+O(M))/O(M),$ and $A(M)/((W+O(M)/O(M))$
is isomorphic to $A(W/M)$ as $A(V/I)$-bimodules.
\end{lem}

\bp{p4.7} Let $j=n-1-(k-1)t$ be an admissible weight. Then
the $A(L(\ell,0))$-bimodule $A(L(\ell,j))$ is
isomorphic to the quotient space of ${\C}[x,y]$ modulo the subspace
$${\C}[x,y]y^{n}+{\C}[x]f_{j,0}(x,y)+{\C}[x]f_{j,1}(x,y)+\cdots
+{\C}[x]f_{j,n-1}(x,y)$$
where
$\displaystyle{f_{j,i}(x,y)=y^{i}\prod_{r=0}^{p-n-1}\prod_{s=0}^{q-k}(x-r-i+st)}.$ The left and right actions of $A(L(\ell,0))$ on $A(L(\ell,j))$ are
given by (\ref{d4.11}).

\ep

{\bf Proof.} First, $M(\ell,j)\simeq U(N_{-})$ as a vector space. Recall that
$B_{0}={\C}(x^{-1}+1)\otimes f+(x^{-2}+x^{-1}){\C}[x^{-1}]\otimes \frak g$.
Since $C=B_{0}\oplus {\C}x^{-1}\otimes e$, by Proposition \ref{p4.1}
$$O(M(\ell,j))=CM(\ell,j)\simeq B_{0}U(N_{-})+ e(-1)U(N_{-}).$$
Since $B_{0}$ is an ideal of $N_{-}$,
$U(N_{-})B_{0}=B_{0}U(N_{-})$ is an ideal of $U(N_{-})$.
Set $L_{0}=N_{-}/B_{0}$. Recall from Section 2 that
$T_+=e(-1)+B_0$ $T_-=f+B_0$ and $T_0=h(-1)+B_0.$ Then $L_0$ is a Lie algebra
spanned by $T_+,T_-,T_0$ and isomorphic to $\frak g$ (see (\ref{d2.13})).

Recall from Proposition \ref{pmff} that $v_{j,1}, v_{j,2}$ are the two
singular vectors of $M(\ell,j)$. Then by Lemma \ref{l4.4} and
Proposition \ref{p4.1}
we have
\begin{eqnarray}
& &A(L(\ell,j))\simeq
M(\ell,j)/(CM(\ell,j)+U(N_{-})v_{j,1}+U(N_{-})v_{j,2})\nonumber\\
& & \ \simeq
U(N_{-})/(B_{0}U(N_{-})+e(-1)U(N_{-})+U(N_{-})F_{1}(n,k)+U(N_{-})F_{2}(n,k))
\end{eqnarray}
as $A(L(\ell,0))$-bimodules. Note that $U(N_-)/B_0U(N_-)\cong U(L_0).$ Thus
\begin{eqnarray}
A(L(\ell,j))\simeq U(L_{0})/(U(L_{0})P(F_{1}(n,k))+U(L_{0})P(F_{2}(n,k))
+T_{+}U(L_{0})).
\end{eqnarray}

For any nonnegative integers $a,b,d$, using Proposition \ref{pf3},
(\ref{d2.14}) and the fact that $G_{\a}=T_+T_--(\a+1)T_0+\a(\a+1)$ we obtain
\begin{eqnarray}
& &T_{+}^{a}T_{0}^{b}T_{-}^{d}P(F_{1}(n,k))\nonumber\\
&=&T_{+}^{a}T_{0}^{b}T_{-}^{d}\left(\prod_{r=0}^{n-1}\prod_{s=1}^{k-1}G_{r+st}\right)T_{-}^{n}\nonumber\\
&=&T_{+}^{a}\left(\prod_{r=0}^{n-1}\prod_{s=1}^{k-1}G_{r+d+st}\right)T_{0}^{b}T_{-}^{d+n}\nonumber\\
&=&T_{+}^{a}\left(\prod_{r=0}^{n-1}\prod_{s=1}^{k-1}(T_{+}T_{-}-(r+d+1+st)T_{0}+(r+d+st)(r+d+1+st)\right)T_{0}^{b}
T_{-}^{n+d}\nonumber\\
&\equiv&T_{+}^{a}\left(\prod_{r=0}^{n-1}\prod_{s=1}^{k-1}(-r-d-1-st)(T_{0}-r-d-st)\right)T_{0}^{b}
T_{-}^{n+d}\;\;\;\mbox{mod }T_{+}U(L_{0}).
\end{eqnarray}
Noticing that $-r-d-1-st\ne 0$ for any $0\le r\le n-1,1\le s\le k-1, d\in
{\Z}_{+}$ we obtain
\begin{eqnarray}
& &U(L_{0})P(F_{1}(n,k))+T_{+}U(L_{0})\nonumber\\
&=&T_{+}U(L_{0})+\sum_{d=0}^{\infty}
{\C}[T_{0}]\left(\prod_{r=0}^{n-1}\prod_{s=1}^{k-1}(T_{0}-r-d-st)\right)T_{-}^{n+d}.
\end{eqnarray}
Similarly, let  $a,b,d$ be any nonnegative integers. If $d<p-n$, we have
\begin{eqnarray}
& &T_{+}^{a}T_{0}^{b}T_{-}^{d}P(F_{2}(n,k))\nonumber\\
&=&T_{+}^{a}T_{0}^{b}T_{-}^{d}T_{+}^{p-n}\prod_{r=1}^{p-n}\prod_{s=1}^{q-k}G_{p-n-r-st}\nonumber\\
&=&T_{+}^{a}T_{0}^{b}\left(\prod_{i=0}^{d-1}
G_{i}\right)T_{+}^{p-n-d}\prod_{r=1}^{p-n}\prod_{s=1}^{q-k}G_{p-n-r-st}\nonumber\\
&=&T_{+}^{a}T_{0}^{b}T_{+}^{p-n-d}\prod_{r=1}^{p-n}\prod_{s=1}^{q-k}\prod_{i=0}^{d-1}G_{p-n-r-st}G_{i+p-n-d}
\nonumber\\
&=&T_{+}^{a+p-n-d}(T_{0}-2(a+p-n-d))^{b}\prod_{r=1}^{p-n}\prod_{s=1}^{q-k}\prod_{i=0}^{d-1}G_{p-n-r-st}G_{i+p-n-d}
\nonumber\\
&\equiv&0\;\;\;\mbox{mod }T_{+}U(L_{0}).
\end{eqnarray}
If $d=m+ p-n$ for some $m\in {\Z}_{+}$, we have
\begin{eqnarray}
& &T_{+}^{a}T_{0}^{b}T_{-}^{d}P(F_{2}(n,k))\nonumber\\
&=&T_{+}^{a}T_{0}^{b}T_{-}^{m}\prod_{i=0}^{p-n-1}
G_{i}\prod_{r=1}^{p-n}\prod_{s=1}^{q-k}G_{p-n-r-st}\nonumber\\
&=&T_{+}^{a}T_{0}^{b}T_{-}^{m}\prod_{r=1}^{p-n}\prod_{s=0}^{q-k}G_{p-n-r-st}\nonumber\\
&=&T_{+}^{a}\left(\prod_{r=1}^{p-n}\prod_{s=0}^{q-k}G_{p+m-n-r-st}\right)T_{0}^{b}T_{-}^{m}\nonumber\\
&=&T_{+}^{a}\left(\prod_{r=0}^{p-n-1}\prod_{s=0}^{q-k}G_{m+r-st}\right)T_{0}^{b}T_{-}^{m}\nonumber\\
&\equiv&T_{+}^{a}\left(\prod_{r=0}^{p-n-1}\prod_{s=0}^{q-k}(-m-r-1+st)(T_{0}-m-r+st)
\right)T_{0}^{b}T_{-}^{m}\;\;\;\mbox{mod }T_{+}U(L_{0}).
\end{eqnarray}
Since $-r-m-1+st\ne 0$ for any $0\le r\le p-n-1,0\le s\le q-k$, we obtain
\begin{eqnarray}
& &U(L_{0})P(F_{2}(n,k))+T_{+}U(L_{0})\nonumber\\
&=&T_{+}U(L_{0})+\sum_{m=0}^{\infty}
{\C}[T_{0}]\left(\prod_{r=0}^{p-n-1}\prod_{s=0}^{q-k}(T_{0}-m-r+st)\right)T_{-}^{m}.
\end{eqnarray}
Thus
\begin{eqnarray*}
& &U(L_{0})P(F_{1}(n,k))+U(L_{0})P(F_{2}(n,k))+T_{+}U(L_{0})\nonumber\\
&\subset &T_{+}U(L_{0})+U(L_{0})T_{-}^{n}+\sum_{i=0}^{n-1}
{\C}[T_{0}]\left(\prod_{r=0}^{p-n-1}\prod_{s=0}^{q-k}(T_{0}-i-r+st)\right)T_{-}^{i}.
\end{eqnarray*}

On the other hand, since $r+d+st\ne m+r'-s't$ for any $0\le r\le n-1, 1\le s
\le k-1, 0\le r'\le p-n-1, 0\le s'\le q-k,d,m\in {\Z}_{+}$,
$\prod_{r=0}^{n-1}\prod_{s=1}^{k-1}(x-r-d-st)$ and
$\prod_{r=0}^{p-n-1}\prod_{s=0}^{q-k}(x-m-r+st)$ are relatively prime. Then
we obtain
$${\C}[T_{0}]T_{-}^{n+i}\subseteq
U(L_{0})P(F_{1}(n,k))+U(L_{0})P(F_{2}(n,k))+T_{+}U(L_{0})$$
for any $i\in {\Z}_{+}$.
This shows that
\begin{eqnarray*}
& &U(L_{0})P(F_{1}(n,k))+U(L_{0})P(F_{2}(n,k))+T_{+}U(L_{0})\nonumber\\
&\supset &T_{+}U(L_{0})+U(L_{0})T_{-}^{n}+\sum_{i=0}^{n-1}
{\C}[T_{0}]\left(\prod_{r=0}^{p-n-1}\prod_{s=0}^{q-k}(T_{0}-i-r+st)\right)T_{-}^{i}.
\end{eqnarray*}
Set $x=T_{0}, y=T_{-}$. Then the proposition follows from Proposition
\ref{p4.3} and Lemma \ref{l4.4}.
$\;\;\;\;\Box$

\bt{t4.8} For admissible weights $j_{i}=n_{i}-1-(k_{i}-1)t$ $(i=1,2)$, the
fusion rules
are given as
follows:
\begin{eqnarray}\label{ef}
L(\ell,j_{1})\times L(\ell,j_{2})
=\sum_{i={\rm max}\{0,n_{1}+n_{2}-p\}}^{{\rm
min}\{n_{1}-1,n_{2}-1\}}L(\ell,j_{1}+j_{2}-2i)
\end{eqnarray}
if $0\le k_{2}-1\le q-k_{1}$, and $L(\ell,j_{1})\times L(\ell,j_{2})=0$
otherwise.
\et

{\bf Proof.} For any admissible weight $j$, let ${\C}v_{j}$ be the
one-dimensional module for Lie
algebra ${\C}h$ such that $hv_{j}=jv_{j}$. Then ${\C}v_{j}$ is the lowest
weight space of $L(\ell,j)$.
By Theorem 3.4 we need to calculate the $A(L(\ell,0))$-module
$A(L(\ell,j_{1}))\otimes _{A(L(\ell,0))}{\C}v_{j_{2}}$.
Using Proposition \ref{p4.7} we get
$$A(L(\ell,j_{1}))\otimes _{A(L(\ell,0))}{\C}v_{j_{2}}\simeq {\C}[x,y]/J$$
where $J$ is the subspace of $\C[x,y]$ spanned by
\begin{equation}\label{d4.24}
\{ x-j_{2},\C[x,y]y^{n_{1}}, f_{j_{1},i}(j_{2},1)\C[x]y^i, i=0,1,
\cdots,n_{1}-1\}
\end{equation}
If $j_{2}$ does not satisfy the relation
$0\le k_{2}-1\le q-k_{1}$, then
$$f_{j_{1},i}(j_{2},1)=\prod_{r=0}^{p-n_{1}-1}\prod_{s=0}^{q-k_{1}}(j_{2}-r-i+st)\ne 0$$
for $0\le i\le n_{1}-1$. Thus $A(L(\ell,j_{1}))\otimes
_{A(L(\ell,0))}{\C}v_{j_{2}}=0$ so that all
the corresponding fusion rules are zero.

Suppose $0\le k_{2}-1\le q-k_{1}$. As before $\C[x]y^i=0$ in $\C[x,y]/J$
if $f_{j_{1},i}(j_{2},1)\ne 0.$ Notice that
$\displaystyle{f_{j_{1},i}(j_{2},1)=\prod_{r=0}^{p-n_{1}-1}\prod_{s=0}^{q-k_{1}}(j_{2}-r-i+st)= 0}$
if and only if $j_{2}-r-i+st= 0$ for some $0\le r\le p-n_{1}-1,0\le s\le
q-k_{1}$. This implies that $0\leq r+i\leq p-2.$ It follows from Remark
\ref{r2.7}
that $r+i=n_2-1.$
That is, $n_{1}+n_{2}-p\le i\le n_{2}-1$.
Therefore
$${\rm max}\{0,n_{1}+n_{2}-p\}\le i\le {\rm min}\{n_{1}-1,n_{2}-1\}.$$
If $n_{1}+n_{2}-p\le i\le n_{2}-1$, then $\C[x]y^{i}$ is not zero in
$\C[x,y]/J.$
Thus
$$\C[x,y]/J\cong \oplus_{{\rm max}\{0,n_{1}+n_{2}-p\}\le i\le{\rm
min}\{n_{1}-1,n_{2}-1\}}\C y^i.$$
{}From (\ref{d4.11}) we get  $x* y^{i}=(j_{2}+j_{1}-2i)y^{i},$ as required.
$\;\;\;\;\Box$

\begin{rem}\label{r4.9} (a) Since $L_z(-1)=L(-1)$ the fusion rules among the
admissible
modules with respect two different operator algebra structure of $L(\ell,0)$
are the same. Thus the fusion rules obatined in Theorem
\ref{t4.8} are also those with respect to the old vertex operator
algebra structure.

(b) After changing the notations one immediately sees that our results
agree with Bernard and Felder's results [BF] on fusion rules by using
BRST cohomology.

(c) Suppose that $\ell$ is an integer. That is, $q=1$
and $p=\ell+2$. Since $j_{i}=n_{i}-1$, we have
$n_{1}+n_{2}-p=j_{1}+j_{2}-\ell$. Since $k_{i}=1$ for any $i$, $0\le
k_{2}-1\le q-k_{1}$ holds automatically. Then
\begin{eqnarray}
L(\ell,j_{1})\times L(\ell,j_{2})
&=&\sum_{i={\rm max}\{0,n_{1}+n_{2}-p\}}^{{\rm
min}\{n_{1}-1,n_{2}-1\}}L(\ell,j_{1}+j_{2}-2i)
\nonumber\\
&=&\sum_{i=0, i\ge j_{1}+j_{2}-\ell}^{{\rm
min}\{j_{1},j_{2}\}}L(\ell,j_{1}+j_{2}-2i)
\nonumber\\
&=&\sum_{j=|j_{1}-j_{2}|, j+j_{1}+j_{2}\le 2\ell}^{j_{1}+j_{2}}L(\ell,j).
\end{eqnarray}
This is exactly the well-known fusion formula (cf. [GW], [TK]).
\end{rem}

\begin{prop}\label{p4.9} Let $M$ be any $V(\ell,\C)$-module. Then
$$C_{2}(M)=({\C}x^{-1}\otimes e+{\C}x^{-1}\otimes f+x^{-2}{\C}[x^{-1}]\otimes
\frak g)M.$$
\end{prop}

{\bf Proof.} Since ${\wt }h=1,{\rm wt}e, {\rm wt}f\notin {\Z}$, by the
definition of $C_{2}(M)$
we get
\begin{eqnarray}
({\C}x^{-1}\otimes e+{\C}x^{-1}\otimes f+x^{-2}{\C}[x^{-1}]\otimes \frak
g)M\subseteq C_{2}(M).
\end{eqnarray}
Set $B_{1}={\C}x^{-1}\otimes e+{\C}x^{-1}\otimes f+x^{-2}{\C}[x^{-1}]\otimes
\frak g$. Let $a$
be a homogeneous element of $V(\ell,{\C})$ such that
\begin{eqnarray}\label{d4.26}
{\rm Res}_{z_{2}}z_{2}^{-n-\varepsilon(a)}Y(a,z_{2})M\subseteq B_{1}M
\end{eqnarray}
for any positive integer $n$. Then for any $k, n\in {\N},u\in M, b\in
\{e,f,h\}$, we have
\begin{eqnarray}
& &{\rm Res}_{z_{2}}z_{2}^{-n-\varepsilon (b(-k)a)}
Y(b(-k)a,z_{2})u\nonumber\\
&=&{\rm Res}_{z_{0}}{\rm Res}_{z_{2}}
z_{2}^{-n-\varepsilon (b(-k)a)}z_{0}^{-k}Y(Y(b,z_{0})a,z_{2})u\nonumber\\
&=&{\rm Res}_{z_{1}}{\rm Res}_{z_{2}}
z_{2}^{-n-\varepsilon(b(-k)a)}(z_{1}-z_{2})^{-k}Y(b,z_{1})
Y(a,z_{2})u\nonumber\\
& &-{\rm Res}_{z_{1}}{\rm Res}_{z_{2}}
z_{2}^{-n-\varepsilon (b(-k)a)}(-z_{2}+z_{1})^{-k}Y(a,z_{2})
Y(b,z_{1})u\nonumber\\
&=&{\rm Res}_{z_{2}}\sum_{i=0}^{\infty}{-k\choose i}
(-z_{2})^{i}z_{2}^{-n-\varepsilon (b(-k)a)}
b(-k-i)Y(a,z_{2})u\nonumber\\
& &-{\rm Res}_{z_{2}}\sum_{i=0}^{\infty}{-k\choose i}
(-1)^{k+i}z_{2}^{-n-k-i-\varepsilon(b(-k)a)}
Y(a,z_{2})b(i)u\nonumber\\
&\equiv&{\rm Res}_{z_{2}}z_{2}^{-n-\varepsilon
(b(-k)a)}b(-k)Y(a,z_{2})u\;\;\mbox{ mod }W\nonumber\\
&\equiv&0\;\;\mbox{ mod }B_1W.\label{d4.27}
\end{eqnarray}
Clearly (\ref{d4.26}) holds for $a={\bf 1}.$ Note that
$V(\ell,0)=U(x^{-1}\C[x]\otimes \frak g){\bf 1}.$ It follows from (\ref{d4.27})
that (\ref{d4.26}) holds for all $a\in V(\ell,\C).$ The proof is complete.
$\;\;\;\;\Box$

\bt{t4.10} The commutative associative algebra
$A_{2}(L(\ell,0),\omega_{z})$ is isomorphic to the quotient algebra
${\C}[x]/\<x^{(p-1)q}\>$.
Consequently, $(L(\ell,0),\omega_{z})$ is $C_{2}$-finite.
\et

{\bf Proof.} First, notice that the Verma module $M(\ell,0)$ is
linearly isomorphic to $U(N_-).$ Recall from Section 2 that
$B_{2}={\C}x^{-1}\otimes f+x^{-2}{\C}[x^{-1}]\otimes \frak g$ is an
ideal of $N_{-},$ $L_{2}=N_{-}/B_{2},$ is the corresponding quotient
Lie algebra spanned by $\bar e=e(-1)+B_2, \bar{f}=f(0)+B_2, \bar{h}=h(-1)+B_2$
and with the commutation relations
\begin{equation}\label{d4.28}
[\bar{e}(-1),\bar{f}(0)]=\bar{h}(-1), [\bar{h}(-1),
\bar{e}(-1)]=[\bar{h}(-1), \bar{f}(0)]=0.
\end{equation}
By
Proposition \ref{p4.9}, we get
$$C_2(M(\ell,0))=B_{2}M(\ell,0)+
e(-1)M(\ell,0)\simeq B_{2}U(N_{-})+ e(-1)U(N_{-}).$$ One easily verifies that
$$A_{2}(L(\ell,0))\simeq U(L_{2})/(
\bar{e}U(L_{2})+U(L_{2})\bar{f}+U(L_{2})P_{2}(F_{2}(1,1))).$$

For any $a,b,m\in {\Z}_{+}$ and $m\ge p-1$ we obtain from Proposition
\ref{pf2} that
\begin{eqnarray}
& &\bar{e}^{a}\bar{h}^{b}\bar{f}^{m}P_{2}(F_{2}(1,1))\nonumber\\
&=&\bar{e}^{a}\bar{h}^{b}\bar{f}^{m}\bar{e}^{p-1}\prod_{r=1}^{p-1}\prod_{s=1}^{q-1}\bar{H}_{p-1-r-st}
\nonumber\\
&=&\bar{e}^{a}\bar{h}^{b}\bar{f}^{m-p+1}\prod_{i=0}^{p-2}\bar{H}_{i}\prod_{r=1}^{p-1}
\prod_{s=1}^{q-1}\bar{H}_{p-1-r-st}
\nonumber\\
&=&\bar{e}^{a}\bar{h}^{b}\left(\prod_{i=0}^{p-2}\prod_{r=1}^{p-1}
\prod_{s=1}^{q-1}\bar{H}_{m-r-st}\bar{H}_{m-p+1+i}\right)
\bar{f}^{m}\nonumber\\
&=&\bar{e}^{a}\bar{h}^{b}\left(\prod_{r=0}^{p-2}\prod_{s=0}^{q-1}\bar{H}_{r+m-p+1-st}\right)\bar{f}^{m-p+1}.
\end{eqnarray}
Here we used the relations $\bar H_{\a}\bar e=\bar e\bar H_{\a+1},$
$\bar H_{\a}\bar f=\bar f\bar H_{\a+1}$ and $\bar f^s\bar e^s=\bar H_0\cdots
\bar H_{s-1}$ which follows from the definition of $\bar H_{\a}$ and the
commutator relations (\ref{d4.28}).
Thus if $a>0$ or $m>p-1$, then
$\bar{e}^{a}\bar{h}^{b}\bar{f}^{p-1+m}P_{2}(F_{2}(1,1))\in
\bar{e}U(L_{2})+U(L_{2})\bar{f}$.

If $a=0$ and $m=p-1$
we have
\begin{eqnarray}
& &\bar h^b\bar{f}^{p-1}P_{2}(F_{2}(1,1))\nonumber\\
&=&\bar h^b\prod_{r=0}^{p-2}\prod_{s=0}^{q-1}\bar{H}_{r-st}\nonumber\\
&\equiv&\bar
h^b\prod_{r=0}^{p-2}\prod_{s=0}^{q-1}(-r-1+st)\bar{h}\;\;\;\mbox{mod }
\bar{e}U(L_{2}).
\end{eqnarray}

Similarly, if $m<p-1$, for any $a,b\in {\Z}_{+}$ we get
$$\bar{e}^{a}\bar{h}^{b}\bar{f}^{m}P_{2}(F_{2}(1,1))\in \bar{e}U(L_{2}).$$

Since $-r-1+st\ne 0$ for any $1\le r\le p-1, 0\le s\le q-1$, we obtain
$$\bar{e}U(L_{2})+U(L_{2})\bar{f}+U(L_{2})P_{2}(F_{2}(1,1))=\bar{e}U(L_{2})+U(\bar{h})\bar{f}+U(L_{2})
\bar{h}^{(p-1)q}.$$
Then the theorem follows if we set $x=\bar{h}$.$\;\;\;\;\Box$

\section{Modular invariance property}
In this section we study modular invariance property of the space linearly
spanned by all
characters
$tr_{L(\ell,j)}e^{2\pi i\tau (L(0)-{1\over 2}zh(0)-{1\over
24}(\frac{3\ell}{\ell+2}-6\ell z^{2}))}$,
where im$\tau>0, z\in {\Q}, 0<z<1$.
In this section we shall first find a modular transformation formula for the
modified characters
for admissible modules.

Following [K] or [KW1]-[KW2], for $m,n\in {\Z}, m>0$ we define
\begin{eqnarray}
\theta_{n,m}(\tau,z)=\sum_{j\in {\Z}+{n\over 2m}}e^{2m\pi i\tau
(j^{2}+jz)},\;\;\; z\in {\C}.
\end{eqnarray}
Set
\begin{eqnarray}
\Theta_{n,m}(\tau)=\sum_{j\in {\Z}+{n\over 2m}}e^{2m\pi i\tau j^{2}}.
\end{eqnarray}
Then
\begin{eqnarray}
\theta_{n,m}(\tau,z)&=&e^{2m\pi i\tau (-{1\over 4}z^{2})}\sum_{j\in
{\Z}+{n\over 2m}}
e^{2m\pi i\tau (j+{1\over 2}z)^{2}}\nonumber\\
&=&e^{-{1\over 2}mz^{2}\pi i\tau}\sum_{j\in {\Z}+\frac{n+mz}{2m}}
e^{2m\pi i\tau j^{2}}.
\end{eqnarray}
Suppose that $z={v\over u}$ is a rational number with $u>0$. Then
\begin{eqnarray}
\theta_{n,m}(\tau,z)&=&e^{-{1\over 2}mz^{2}\pi i\tau}\sum_{j\in
{\Z}+\frac{nu+mv}{2mu}}
e^{2m\pi i\tau j^{2}}\nonumber\\
&=&e^{-{1\over 2}mz^{2}\pi i\tau}\Theta_{nu+mv,mu}({\tau\over u}).
\end{eqnarray}

As in Section 4, we let $\ell=-2+{p\over q}$ be a fixed admissible level, where
$p\ge 2 ,q$ are
relatively prime positive integers.
Let $P_{\ell}$ be the set of all admissible weights (mod ${\C}\delta$) of level
$\ell$. Then
$$P_{\ell}=\{ j=n-kt| n,k\in {\Z}_{+}, n\le p-2, k\le q-1\}.$$
Set
$c_{\ell}=\frac{3\ell}{\ell+2}$. For any rational number $z$, we set
$c_{\ell,z}=c_{\ell}-6\ell z^{2}$.
In Section 4 we have studied the vertex operator algebra or chiral algebra
$L(\ell,0)$ under a different Virasoro vector $\omega_{z}$ which
has a central charge $c_{\ell,z}$. That is, the rank of
$(L(\ell,0),\omega_{z})$ is
$c_{\ell,z}$. With this motivation we define the following characters
\begin{eqnarray}
\chi_{j}(\tau,z):=\tr_{L(\ell,j)}e^{2\pi i\tau (L_{z}(0)-{1\over
24}c_{\ell,z})}
=\tr_{L(\ell,j)}e^{2\pi i\tau (L(0)-{1\over 2}zh(0)-{1\over 24}c_{\ell,z})}.
\end{eqnarray}
For an admissible weight $j=n-kt\in P_{\ell}$, set
$$a=pq, b_{j}^{\pm}=q(\pm (n+1)-kt).$$
Now we restrict $z$ to be a positive rational number less than $1$.

\begin{rem}\label{r5.1}
In [KW1]-[KW2], the following defined character has been considered:
\begin{eqnarray}
\bar{\chi}_{j}(\tau,z)=\tr_{L(\ell,j)}e^{2\pi i\tau(L(0)-{1\over
2}zh(0)-{1\over 24}c_{\ell})},
\end{eqnarray}
and it was proved that
\begin{eqnarray}
\bar{\chi}_{j}(\tau,z)=
\frac{\theta_{b_{j}^{+},a}(\tau,q^{-1}z)-\theta_{b_{j}^{-},a}(\tau,q^{-1}z)}
{\theta_{1,2}(\tau,z)-\theta_{-1,2}(\tau,z)}.
\end{eqnarray}
\end{rem}

Using KW's character formula we obtain
\begin{eqnarray}
& &\ \ \ \chi_{j}(\tau,z)\nonumber\\
& &=e^{{1\over 2}\ell z^{2}\pi i\tau }\bar{\chi}_{j}(\tau,z)\nonumber\\
& &=e^{{1\over 2}\ell z^{2}\pi i\tau }
\frac{\theta_{b_{j}^{+},a}(\tau,q^{-1}z)-\theta_{b_{j}^{-},a}(\tau,q^{-1}z)}
{\theta_{1,2}(\tau,z)-\theta_{-1,2}(\tau,z)}\nonumber\\
& &=e^{{1\over 2}\ell z^{2}\pi i\tau }
e^{-{1\over 2}aq^{-2}z^{2}\pi i\tau}e^{z^{2}\pi i\tau}
\frac{\Theta_{qub_{j}^{+}+av,aqu}({\tau\over
qu})-\Theta_{b_{j}^{-}qu+av,aqu}({\tau\over qu})}
{\Theta_{u+2v,2u}({\tau\over u})-\Theta_{-u+2v,2u}({\tau\over u})}\nonumber\\
& &=e^{{1\over 2}z^{2}\pi i\tau (\ell+2-aq^{-2})}
\frac{\Theta_{qub_{j}^{+}+av,aqu}({\tau\over
qu})-\Theta_{uqb_{j}^{-}+av,aqu}({\tau\over qu})}
{\Theta_{u+2v,2u}({\tau\over u})-\Theta_{-u+2v,2u}({\tau\over u})}\nonumber\\
& &=\frac{\Theta_{uqb_{j}^{+}+av,aqu}({\tau\over
qu})-\Theta_{qub_{j}^{-}+av,aqu}({\tau\over qu})}
{\Theta_{u+2v,2u}({\tau\over u})-\Theta_{-u+2v,2u}({\tau\over u})}.
\end{eqnarray}
Then $\ch_{j}$ is a modular function with $c_{\ell,z}$ as the modular anomaly
rather than $c_{\ell}$.

\begin{rem}\label{r5.2} In [KW1] the following transformation law was given:
\begin{eqnarray}\label{e5.3}
\bar{\chi}_{j}(-\tau^{-1},\tau z)={1\over 2i}\sqrt{{2\over a}}
\sum_{j'\in P_{\ell}} \left(e^{-i\pi b_{+}b_{-}'/a}-e^{-i\pi
b_{+}b_{+}'/a}\right)
\bar{\chi}_{j'}(\tau,z).
\end{eqnarray}
Later in [KW2], a correction was made
by adding the factor $e^{{1\over 2}\ell z^{2}\pi i\tau}$ on the right-hand side
of (\ref{e5.3}). That is,
\begin{eqnarray}\label{e5.4}
\bar{\chi}_{j}(-\tau^{-1},\tau z)={1\over 2i}\sqrt{{2\over a}}e^{{1\over 2}\ell
z^{2}\pi i\tau}
\sum_{j'\in P_{\ell}} \left(e^{-i\pi b_{+}b_{-}'/a}-e^{-i\pi
b_{+}b_{+}'/a}\right)
\bar{\chi}_{j'}(\tau,z).
\end{eqnarray}
Based on the modular transformation law ((\ref{e5.3}) without the factor
$e^{{1\over 2}\ell z^{2}\pi i\tau}\!$), the fusion
rules have been calculated in [KS] and [MW] by using Verlinde formula
[V]. Unfortunately, some of them are negative. On the other hand, the
correct formula (\ref{e5.4}) can not be used to compute the fusion
because the coefficients in (\ref{e5.4}) involve the variable $\tau.$
This puzzles both
mathematicians and physicists.
\end{rem}

For a $\Z$-graded rational vertex operator algebra satisfying $C_2$
condition and the Virasoro condition it is proved in [Z] that
the space spanned by $\tr_Mq^{L(0)-{c\over 24}}$ for all
irreducible modules $M$ modular invariant. If the Virasoro condition
is replaced by the primary field condition (cf. Remark \ref{r3.12})
one still has the modular invariance of the space by
 modifying Zhu's proof.
Now we have a $\Q$-graded
rational vertex operator algebra $(L(\ell,0),\o_z)$ satisfying
$C_2$ condition (see Theorem \ref{t4.10}) and primary field condition
(see Remark \ref{r3.12}). Unfortunately Zhu's modular invariance theorem [Z]
does not apply to $Q$-graded vertex operator algebra.
This raise a question:
Is the space linearly spanned by $\chi_{j}(\tau,z)$ modular invariant
under the transformation $\tau\mapsto -\tau^{-1}$ with $z$ being fixed?
This question will be discussed in our coming paper [DLiM3].

\end{document}